\documentclass[11pt]{amsart}
\usepackage{graphicx}
\usepackage{version}

\excludeversion{longversion}

\numberwithin{equation}{section}

\newcommand{\transpose}[1]{{{#1}^t}}

\newcommand{\abs}[1]{{\lvert#1\rvert}}
\newcommand{\norm}[1]{\lVert#1\rVert}

\newcommand{\half}{\tfrac{1}{2}}

\newcommand{\fourth}{\tfrac{1}{4}}

\newcommand{\ol}[1]{\overline{#1}}

\DeclareMathOperator{\Order}{O}

\DeclareMathOperator{\Span}{span}

\DeclareMathOperator{\adjoint}{adjoint}

\DeclareMathOperator{\cross}{\times}

\DeclareMathOperator{\diag}{diag}

\DeclareMathOperator{\id}{I}

\DeclareMathOperator*{\ord}{ord}

\DeclareMathOperator{\rank}{rank}

\DeclareMathOperator{\suchthat}{|}

\DeclareMathOperator{\tensor}{\otimes}
\DeclareMathOperator{\tr}{tr}
\DeclareMathOperator{\trace}{tr}

\newcommand{\MatrixGroup}[1]{{\rm{#1}}}

\newcommand{\matGL}{\MatrixGroup{GL}}

\newcommand{\matSL}{\MatrixGroup{SL}}

\newcommand{\matSU}{\MatrixGroup{SU}}
\newcommand{\matU}{\MatrixGroup{U}}

\newcommand{\matsl}{\MatrixGroup{sl}}
\newcommand{\matsu}{\MatrixGroup{su}}

\newcommand{\matN}[1]{\MatrixGroup{M}_{#1\times #1}}

\theoremstyle{plain}
\newtheorem{theorem}{Theorem}[section]

\newtheorem{lemma}[theorem]{Lemma}
\newtheorem{proposition}[theorem]{Proposition}
\newtheorem*{corollary*}{Corollary}
\newtheorem*{lemma*}{Lemma}
\newtheorem*{proposition*}{Proposition}
\newtheorem*{theorem*}{Theorem}
\theoremstyle{definition}

\newtheorem{definition}[theorem]{Definition}

\newtheorem{remark}[theorem]{Remark}

\newtheorem*{algorithm*}{Algorithm}
\newtheorem*{application*}{Application}
\newtheorem*{assertion*}{Assertion}
\newtheorem*{assumption*}{Assumption}
\newtheorem*{axiom*}{Axiom}
\newtheorem*{claim*}{Claim}
\newtheorem*{conjecture*}{Conjecture}
\newtheorem*{definition*}{Definition}
\newtheorem*{example*}{Example}
\newtheorem*{notation*}{Notation}
\newtheorem*{note*}{Note}
\newtheorem*{observation*}{Observation}
\newtheorem*{question*}{Question}
\newtheorem*{remark*}{Remark}
\newtheorem*{problem*}{Problem}
\theoremstyle{plain}

\newcommand{\bbC}{\mathbb{C}}

\newcommand{\bbH}{\mathbb{H}}

\newcommand{\bbN}{\mathbb{N}}

\newcommand{\bbP}{\mathbb{P}}

\newcommand{\bbR}{\mathbb{R}}
\newcommand{\bbS}{\mathbb{S}}

\newcommand{\bbZ}{\mathbb{Z}}

\newcommand{\calA}{\mathcal{A}}

\newcommand{\calC}{\mathcal{C}}
\newcommand{\calD}{\mathcal{D}}

\newcommand{\calG}{\mathcal{G}}
\newcommand{\calH}{\mathcal{H}}

\newcommand{\calL}{\mathcal{L}}
\newcommand{\calM}{\mathcal{M}}

\newcommand{\calU}{\mathcal{U}}
\newcommand{\calV}{\mathcal{V}}

\newcommand{\LoopSL}[1]{\Lambda_{#1}\matSL_2(\bbC)}

\newcommand{\LoopuSL}[1]{\Lambda_{#1}^\ast\matSL_2(\bbC)}

\newcommand{\LooppSL}[1]{\Lambda_{#1}^{+}\matSL_2(\bbC)}

\newcommand{\gauge}[2]{{{#1}{.}{#2}}}

\newcommand{\spacecomma}{,}
\newcommand{\spaceperiod}{.}

\newcommand{\Holo}[2]{\calH_{#1}#2}
\newcommand{\Mero}[2]{\calM_{#1}#2}

\numberwithin{equation}{section}

\numberwithin{equation}{subsection}

\newcommand{\Section}[1]{\section{#1}}
\newcommand{\Sectionstar}[1]{\section*{#1}}
\newcommand{\Subsection}[1]{\subsection{#1}}

\newcommand{\Note}[1]{\footnote{#1}}

\newcounter{COUNTER}

\renewcommand{\matGL}[2]{\MatrixGroup{GL}_{#1}{#2}}

\renewcommand{\matSL}[2]{\MatrixGroup{SL}_{#1}{#2}}

\renewcommand{\matSU}[2]{\MatrixGroup{SU}_{#1}{#2}}
\renewcommand{\matU}[2]{\MatrixGroup{U}_{#1}{#2}}

\renewcommand{\matsl}[2]{\MatrixGroup{sl}_{#1}{#2}}
\renewcommand{\matsu}[2]{\MatrixGroup{su}_{#1}{#2}}
\newcommand{\matM}[2]{\MatrixGroup{M}_{#1}{#2}}

\renewcommand{\LoopSL}[1]{\Lambda_{#1}\matSL{2}{\bbC}}
\renewcommand{\LooppSL}[1]{\Lambda_{#1}^{+}\matSL{2}{\bbC}}
\renewcommand{\LoopuSL}[1]{\Lambda_{#1}^\ast\matSL{2}{\bbC}}

\newcommand{\Hom}{\rm{Hom}}

\author{W. Rossman}
\address{Kobe University, Rokko Kobe 657-8501, Japan}
\email{wayne@math.kobe-u.ac.jp}

\author{N. Schmitt}
\address{Institut f\"ur Mathematik,
Technische Universit\"at Berlin,
10623 Berlin, Germany}
\email{nick@gang.umass.edu}

\title%
[Simultaneous unitarizability and $k$-noids]%
{Simultaneous unitarizability of $\matSL{n}{\bbC}$-valued maps,
and constant mean curvature $k$-noid monodromy}

\dedicatory{
\textit{Dedicated to Professor Takeshi Sasaki on his sixtieth birthday}}

\date{\today}

\begin{document}

\typeout{===abstract==============}
\begin{abstract}
We give necessary and sufficient local conditions for the
simultaneous unitarizability of a set of
analytic matrix maps from an analytic 1-manifold
into $\matSL{n}{\bbC}$ under conjugation by a single analytic matrix map.
We apply this result to the monodromy arising from an integrable
partial differential equation
to construct a family of $k$-noids,
genus-zero constant mean curvature surfaces with three or more
ends in Euclidean, spherical and
hyperbolic $3$-spaces.
%
\end{abstract}

\maketitle

\typeout{===intro==============}

\Sectionstar{Introduction}

In this paper
 we find necessary and sufficient conditions
for the existence of an $\matSL{n}{\bbC}$-valued
analytic matrix map on an analytic 1-manifold which
simultaneously unitarizes a given set of analytic matrix maps via conjugation.
We apply these results to construct families of
constant mean curvature (CMC) immersions
with arbitrarily many ends into ambient
3-dimensional space forms
(Theorem~\ref{thm:nnoid}).
The ends are conjectured to be asymptotic to
half-Delaunay surfaces.

We show that the existence of a global unitarizer is equivalent
to the existence of analytic unitarizers
defined only on local neighborhoods (Theorem~\ref{thm:1uni}).
In the case of $\matSL{2}{\bbC}$
the necessary and sufficient conditions for global simultaneous unitarization
are local diagonalizability,
pointwise simultaneous
unitarizability, and pairwise infinitesimal irreducibility
(Theorem~\ref{thm:uni-2by2}).
This latter condition means that at each point of the $1$-manifold, the
coefficient of the leading term of the series expansion of the
commutator has full rank.
For general $\matSL{n}{\bbC}$, global unitarizability is equivalent to
pointwise unitarizability together with a graph condition
(Theorem~\ref{thm:nuni}).
These results are proven by linearizing the unitarization problem
and applying analytic Cholesky decompositions.

The Unitarization Theorem~\ref{thm:nuni}  is a refinement of the variant
$r$-unitarization Theorem~\ref{thm:runi}
(see~\cite{SKKR:spaceforms,DorW:uni} for the case of $\matSL{2}{\bbC}$).
In Theorem~\ref{thm:runi}, the analytic curve is
the standard unit circle $\bbS^1$,
and an analytic simultaneous unitarizer $C$ is found on a
radius-$r$ circle for some $r$ less than $1$.
While the unitarized loops extend holomorphically to $\bbS^1$,
the unitarizer $C$ generally has branch points.
The conditions of the Unitarization Theorem~\ref{thm:nuni} are the
obstructions to extending $C$ holomorphically to $\bbS^1$.

One application of the Unitarization Theorem~\ref{thm:nuni}
is to solving the monodromy
problem arising in the construction of CMC
surfaces via the
extended Weierstrass representation. In this construction, using integral
system methods, the problem of closing the surface is solved by unitarizing
the monodromy group, whose elements are defined on a loop.
The unitarization theorem
provides a construction of a global analytic simultaneous unitarizer once
it is known that the monodromy group is pointwise simultaneously unitarizable
along the loop.

\begin{figure}[ht]
  \centering
  \includegraphics[scale=0.35]{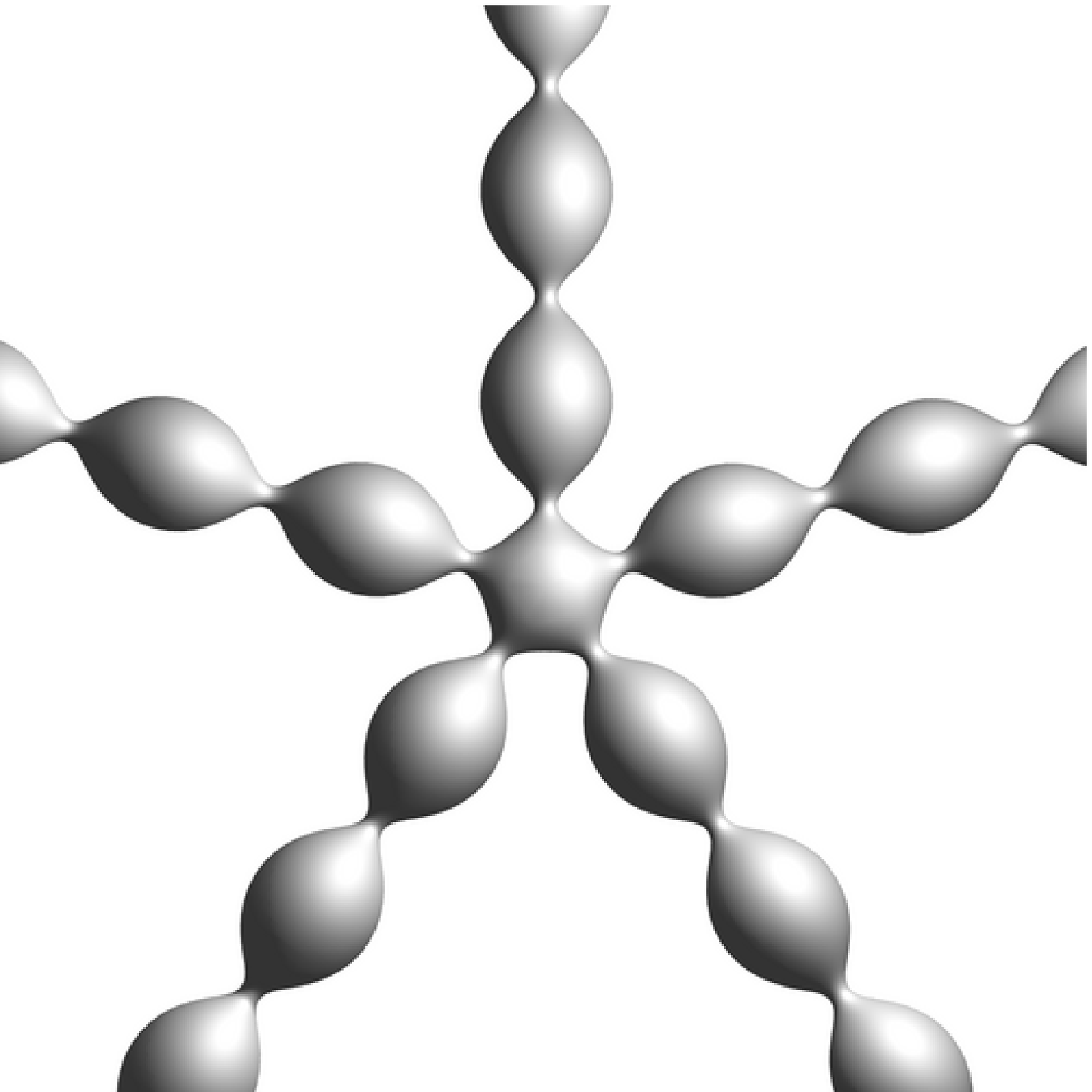}
  \includegraphics[scale=0.35]{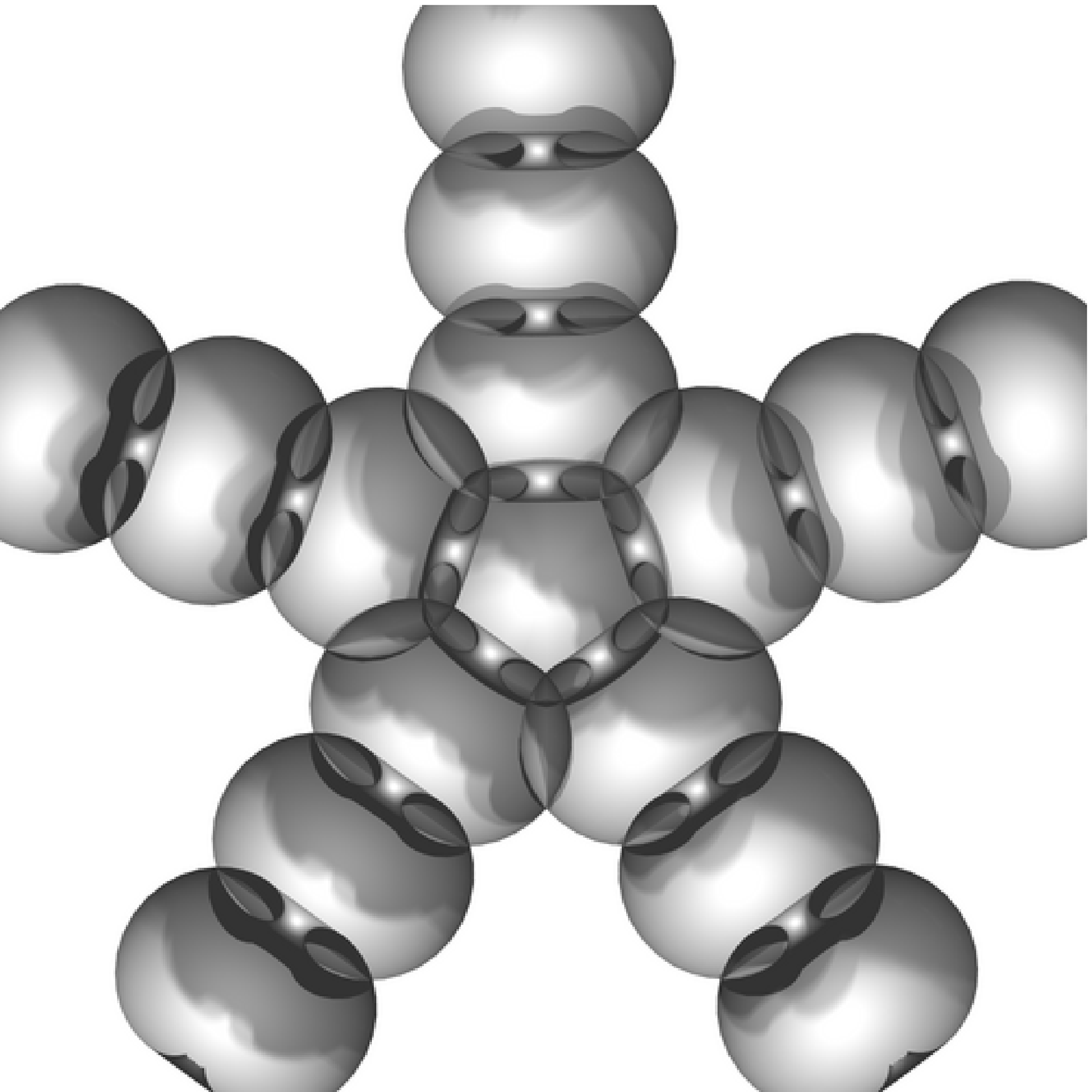}
\caption{
Symmetric CMC $5$-noids in $\bbR^3$,
one with unduloidal ends and one with nodoidal ends (cutaway view).
The images were produced
by \texttt{cmclab}~\cite{Schmitt:cmclab}.
}
\label{fig:nnoid1}
\end{figure}

\begin{figure}[ht]
 \centering
  \includegraphics[scale=0.34]{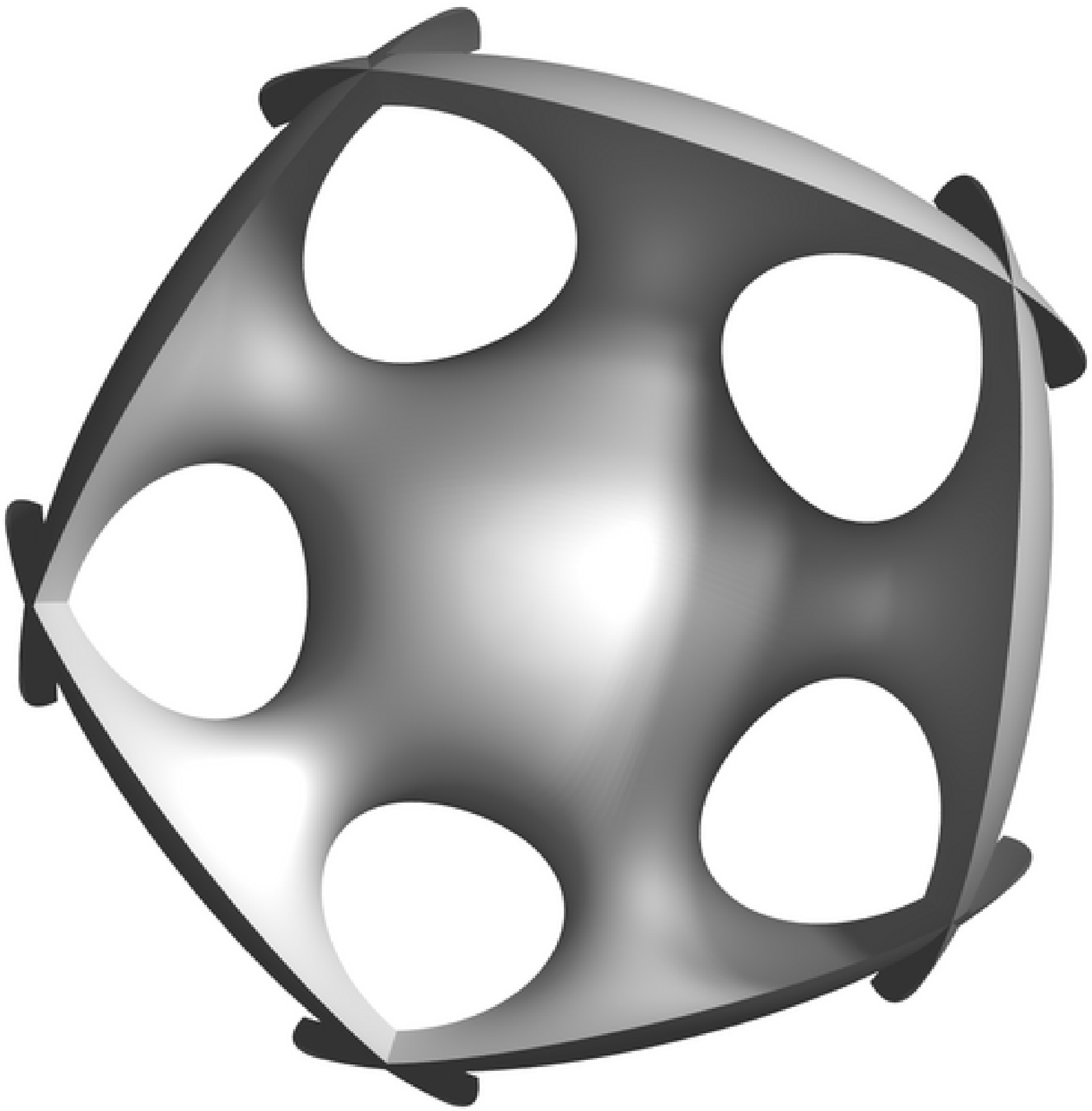}
  \includegraphics[scale=0.36]{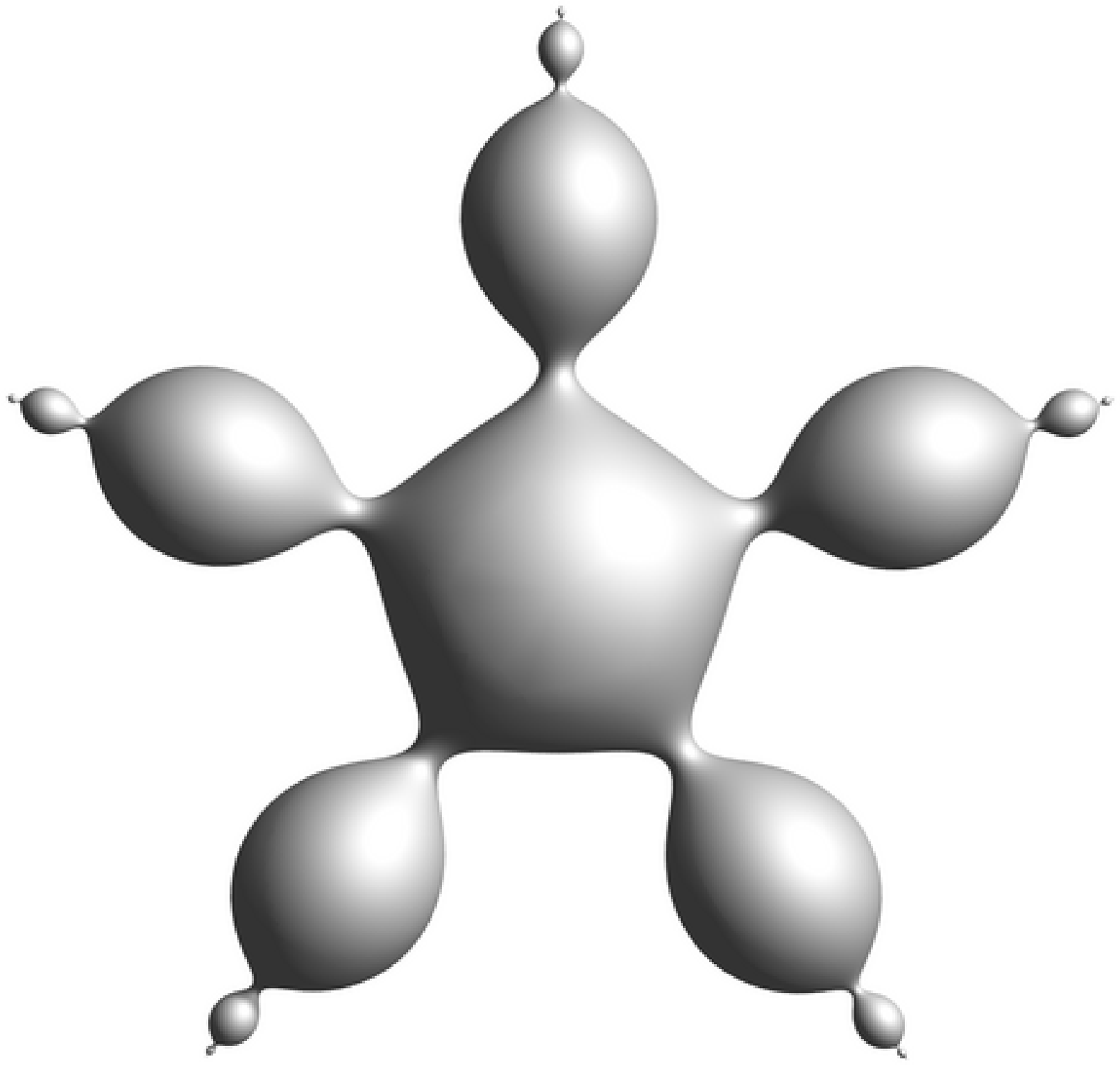}
\caption{
Symmetric CMC $5$-noids in $\bbS^3$ and $\bbH^3$.
Here, $\bbS^3$ has been stereographically projected
to $\bbR^3\cup\{\infty\}$, and $\bbH^3$ is shown
in the Poincar\'{e} model.}
\label{fig:nnoid2}
\end{figure}

Hence in the second part of the paper we apply the unitarization theorem to
the construction of CMC genus-zero
surfaces with arbitrary numbers $k\geq 3$ of ends which are
conjectured to be asymptotic to half-Delaunay surfaces,
lying in ambient $3$-dimensional
space forms (see Figures~\ref{fig:nnoid1} and~\ref{fig:nnoid2}).
We call these surfaces $k$-noids, and trinoids when
$k=3$. For $k=2$ these are the well-known Delaunay surfaces,
CMC surfaces of revolution with
translational periodicity.

Trinoids in $\bbR^3$, either embedded or non-Alexandrov-embedded, were
first constructed by Kapouleas~\cite{Kap1} using techniques that glue parts
of CMC surfaces together via the study of Jacobi operators, and later
there has been further work in this
direction~\cite{Mazzeo-Pacard,Ratzkin:endtoend},
but this approach only gives examples that are in
some sense ``close'' to the boundary of the moduli space of the
surfaces.

A construction that gives a broader collection of
Alexandrov embedded trinoids~\cite{GKS:tri,GKS:coplanar}
and $k$-noids~\cite{GB-newcmc} 
in $\bbR^3$ with embedded ends
asymptotic to Delaunay unduloids
was found by
Gro{\ss}e-Brauckmann, Kusner and Sullivan,
using an isometric correspondence between minimal
surfaces in the $3$-dimensional sphere $\bbS^3$ and CMC $1$
surfaces in $\bbR^3$.
The family of $k$-noids we construct here also include
ends whose potentials are perturbations of Delaunay nodoid potentials.

Constant mean curvature
trinoids in $\bbR^3$ with embedded ends via loop group techniques are
constructed in~\cite{SKKR:spaceforms,DorW:tri}
by methods derived from integrable systems techniques.
Developed initially by
Dorfmeister, Pedit and Wu~\cite{DPW},
%
this construction employs the $r$-Unitarization Theorem
together with the $r$-Iwasawa decomposition~\cite{McI1},
a generalization of the Iwasawa decomposition on the
unit circle to a radius-$r$ circle with $r<1$.

These trinoids extend to larger classes of non-Alexandrov-embedded
trinoids: In one way, Kilian, Sterling and the second
author~\cite{KilSS} dressed these trinoids into ``bubbleton''
versions, also conformal to thrice-punctured spheres;
computer graphics suggest that these bubbleton versions
have ends asymptotic to embedded Delaunay unduloids and
are not Alexandrov embedded.  In another way,
Kilian, Kobayashi and the authors~\cite{SKKR:spaceforms}
extended the class of trinoids to CMC
$1$ surfaces in $\bbR^3$ whose potentials are perturbations of
nonembedded Delaunay nodoid potentials, and hence are conjectured
to be not Alexandrov embedded.  Also, in~\cite{SKKR:spaceforms}, trinoids
in $\bbS^3$ and hyperbolic $3$-space $\bbH^3$ were proved
to exist, including examples that are not Alexandrov embedded.  (The
asymptotic behavior of the ends is considered in a separate work by Kilian
and the authors~\cite{KRS:asymptotics}.)

In this work, we establish the closing conditions for
trinoids and symmetric $k$-noids
via loop group techniques by a more elementary
approach using only the $1$-unitarization
theorem and $1$-Iwasawa decomposition.


\Section{Simultaneous unitarizability for $\matSL{n}{\bbC}$}
\label{sec:uni}

\typeout{===prelim==============}
\Subsection{Preliminaries}

An \emph{analytic curve} is a connected real analytic one-dimensional manifold
without boundary.
%
On an analytic curve $\calC$
we denote by $\Holo{\calC}{\calV}$ the set of
analytic maps $\calC\to\calV$
into a space $\calV$,
and by $\Mero{\calC}{\calV}$ the set of
analytic maps $\calC\to\calV$ with possible poles.

$M \in \Holo{\calC}{\matSL{n}{\bbC}}$ is
\emph{locally diagonalizable at $p\in\calC$}
if there exists a neighborhood $\calU \subseteq \calC$ of $p$
and $V\in\Holo{\calU}{\matSL{n}{\bbC}}$ such that
$VMV^{-1}$ is diagonal.

$M_1,\dots,M_q \in \Holo{\calC}{\matSL{n}{\bbC}}$ are
\emph{locally simultaneously unitarizable at $p\in\calC$} if
there exists a neighborhood $\calU \subseteq \calC$ of $p$
and $V\in\Holo{\calU}{\matSL{n}{\bbC}}$
such that
$VM_1V^{-1},\dots,VM_qV^{-1}\in \Holo{\calU}{\matSU{n}{\bbC}}$.

$M_1,\dots,M_q \in \Holo{\calC}{\matSL{n}{\bbC}}$ are
\emph{simultaneously unitarizable} on $\calC$ if
there exists $V\in\Holo{\calC}{\matSL{n}{\bbC}}$
such that
$VM_1V^{-1},\dots,VM_qV^{-1}\in \Holo{\calC}{\matSU{n}{\bbC}}$.

For $f\in\Mero{\calC}{\bbC}$, define $f^\ast = \ol{f}$,
and for $M\in\Mero{\calC}{\matM{n\times n}{\bbC}}$,
define $M^\ast = \transpose{\ol{M}}$.

For $0<r<s<\infty$, let $\calA_{r,s}\subset\bbC$ denote the open annulus
$\calA_{r,s}=\{\lambda\in\bbC\suchthat r < \abs{\lambda} < s\}$.



A subset $\calG\subset\matM{n\times n}{\bbC}$ is
\emph{reducible} if there exists
a proper non-zero subspace $V\subset\bbC^n$ such that
$GV\subset V$ for all $G\in\calG$.
In the case $\calG=\{A,\,B\}$ is a set of two elements,
we say $A,\,B\in\matM{n\times n}{\bbC}$ are reducible.

We have the following Schur-type lemma.

\begin{lemma}
\label{lem:schur}
If $A,\,B\in\matM{n\times n}{\bbC}$ are irreducible and
$X\in\matM{n\times n}{\bbC}$ commutes with $A$ and $B$,
then $X$ is a multiple of the identity matrix $\id\in\matM{n\times n}{\bbC}$.
\end{lemma}

\begin{proof}
Suppose $X\not\in\bbC\id$ and let $\lambda\in\bbC$
be an eigenvalue of $X$.
Let $V = \{v\in\bbC^n\suchthat Xv = \lambda v\}$,
so $V \ne \{0\}$.
Then $V\subset\bbC^n$ is a proper subspace because $X\not\in\bbC\id$.

For any $v\in V$, we have $XAv = AXv = \lambda A v$, so $Av\in V$.
Hence $AV\subseteq V$. Similarly, $BV\subseteq V$, so
$A$ and $B$ would be reducible.
\end{proof}

\typeout{===linearmap==============}
\Subsection{Linearizing the simultaneous unitarization problem}

We begin with an elementary proof of a
specialization of standard results in the
theory of holomorphic vector bundles,
constructing a global kernel of a
suitable bundle map which depends holomorphically
(or real analytically) on a parameter.

\begin{lemma}
\label{lem:holo-kernel}
Let $\calD$ be an analytic 1-manifold or 2-manifold,
$N,\,M\in\bbN$, and $L:\calD\to\Hom(\bbC^N,\bbC^M)$ a holomorphic map.
Suppose $\dim\ker L=1$ on $\calD$ away from a subset $S\subset\calD$
of isolated points.
Then
\begin{enumerate}
\item
$\dim\ker L\ge 1$ on $\calD$.
\item
There exists a holomorphic map $X:\calD\to\bbC^N$
which is not identically $0$
such that $X\in\ker L$, that is,
for each $p\in\calD$, $X(p)\in\ker L(p)$.
\end{enumerate}
\end{lemma}

\begin{proof}
Since $L$ has rank $N-1$ on $\calD\setminus S$,
all the $N\times N$ minor determinants of $L$
are holomorphic on $\calD$ and zero on $\calD\setminus S$,
and hence $0$ on $\calD$.
Hence $\dim\ker L\ge 1$ on $\calD$.

Since $L$ has rank $N-1$ on $\calD\setminus S$,
then there exists an $(N-1)\times(N-1)$ minor determinant
of $L$ which is not
identically $0$ on $\calD$.
Hence by a permutation we may assume without
loss of generality that
\[
L =
\begin{pmatrix}
A & B \\ C & D
\end{pmatrix},
\]
where $A\in\Holo{\calD}{\matM{(N-1)\times(N-1)}{\bbC}}$
with $a:=\det A\not\equiv 0$
on $\calD$,
$B\in\Holo{\calD}{\matM{(N-1)\times 1}{\bbC}}$,
$C\in\Holo{\calD}{\matM{(M-N+1)\times(N-1)}{\bbC}}$,
and $D\in\Holo{\calD}{\matM{(M-N+1)\times 1}{\bbC}}$.

Define $X=\transpose{(-a A^{-1} B,\,a)}$.
Then $a A^{-1}$ is holomorphic on $\calD$ because
its entries are polynomials in the entries of $A$.
Hence $X$ is holomorphic on $\calD$.
Moreover, $X\not\equiv 0$ because $a\not\equiv 0$.

It is clear that the upper $(N-1)\times 1$ block of $LX$ is $0$.
To show that the lower $(M-N+1)\times 1$ block of $LX$ is zero,
let $A_k\in\Holo{\calC}{\matM{1\times(N-1)}{\bbC}}$
and $B_k\in\Holo{\calC}{\bbC}$ be the respective $k$'th rows of $A$ and $B$
($k\in\{1,\dots,N-1\}$).
Then since $B-AA^{-1}B=0$,
\begin{equation}
\label{eq:kerX2}
B_k - A_k A^{-1}B = 0,\quad
k\in\{1,\dots,N-1\}.
\end{equation}
Fix $j\in\{1,\dots,M-N+1\}$
and let
$C_j\in\Holo{\calC}{\matM{1\times(N-1)}{\bbC}}$
and $D_j\in\Holo{\calC}{\bbC}$
be the respective $j$'th rows of
$C$ and $D$.
Because all $N\times N$ minor determinants of $L$ are $0$, then
$C_j=\sum r_k A_k$
and $D_j = \sum r_k B_k$ for some holomorphic scalar functions
$r_1,\dots,r_{N-1}$.
By~\eqref{eq:kerX2},
$C_j A^{-1}B = D_j$.
Since this holds for all $j\in\{1,\dots,M-N+1\}$, we have
$D - C A^{-1} B = 0$.
Hence the lower $(M-N+1)\times 1$ block
$a(D-CA^{-1}B)$ of $LX$ is $0$,
$LX=0$, and $X\in\ker L$.
\end{proof}

A global simultaneous unitarizer of
a suitable set of analytic maps $M_1,\dots,M_q$ on an analytic curve
is constructed in two steps.
First, a global analytic solution $X$ to the linear system
$XM_kX^{-1}={M_k^\ast}^{-1}$, $k\in\{1,\dots,q\}$, is found
which is Hermitian positive definite.
Then the analytic
Cholesky decomposition gives $X=V^\ast V$,
and $V$ is a simultaneous unitarizer of
$M_1,\dots,M_q$.

We continue by defining a linear map $L$ whose
kernel will contain $X$.

\begin{definition}
\label{def:L}
For $M_1,\dots,M_q\in\matSL{n}{\bbC}$ ($n\ge 1$),
define $L=\calL(M_1,\dots,M_q)$ as
the linear map $L:\matM{n\times n}{\bbC}\to(\matM{n\times n}{\bbC})^n$
\begin{equation*}
L(X) = \left(XM_1-{M_1^\ast}^{-1}X,\dots,XM_q-{M_q^\ast}^{-1}X\right).
\end{equation*}
\end{definition}
\noindent $L$ is similarly defined
for $M_1,\dots,M_q\in\Holo{\calC}{\matSL{n}{\bbC}}$.

The linear map $L(M_1,\dots,M_q)$ has the following easily computed properties.
The first of these properties motivates the definition of $L$.

\begin{lemma}
\label{lem:L}
With $n\ge 2$, let $M_1,\dots,M_q\in\matSL{n}{\bbC}$
and let $L=\calL(M_1,\dots,M_q)$ be as in Definition~\ref{def:L}.
Then
\begin{enumerate}
\item
If $A\in\matSL{n}{\bbC}$ and $A^\ast A\in\ker L$, then
$A$ simultaneously unitarizes $M_1,\dots,M_q$.
\item
If $X\in\ker L$, then $X^\ast\in\ker L$.
\item
Let $C\in\matGL{n}{\bbC}$
and let $\widetilde{L}=\calL(CM_1C^{-1},\dots,CM_qC^{-1})$.
Then $X\in\ker \widetilde{L}$ if and only if 
$C^\ast XC\in\ker L$.
\end{enumerate}
\end{lemma}

We will require the following further properties
of $L$.

\begin{lemma}
\label{lem:dimkerL}
With $n\ge 2$, let $M_1,\dots,M_q\in\matSL{n}{\bbC}$
and $L=\calL(M_1,\dots,M_q)$.
If for some $i,\,j\in\{1,\dots,q\}$, $M_i$ and $M_j$ are irreducible,
and if $M_1,\dots,M_q$ are simultaneously unitarizable, then
$\dim\ker L=1$.
\end{lemma}

\begin{proof}
Let $C$ be a simultaneous unitarizer of $M_1,\dots,M_q$,
define $P_k=CM_kC^{-1}\in\matSU{n}{}$ ($k\in\{1,\dots,n\}$),
and set $\widetilde{L} = \calL(P_1,\dots,P_r)$.
Since $P_k ={P_k^\ast}^{-1}$,
then $\ker \widetilde{L}$ is the set of $X\in\matM{n\times n}{\bbC}$
such that $[X,\,P_k]= 0$ ($k\in\{1,\dots,q\}$).
Since $M_i$ and $M_j$ are irreducible, then
$P_i$ and $P_j$ are irreducible.
By Lemma~\ref{lem:schur},
$\ker\widetilde{L} = \bbC\tensor\id$.
By Lemma~\ref{lem:L}(iii),
$\ker L = \bbC\tensor (C^\ast C)$,
so $\dim\ker L=1$.
\end{proof}

\begin{lemma}
\label{lem:posdef}
Let
$M_1,\,M_2\in\matSL{n}{\bbC}$ and $L=\calL(M_1,\,M_2)$.
Suppose $M_1$ and $M_2$ are irreducible
and simultaneously unitarizable.
Let $X\in\matM{n\times n}{\bbC}\setminus\{0\}$
and suppose $X\in\ker L$ and $X^\ast = X$.
Then $X$ is positive or negative definite.
\end{lemma}

\begin{proof}
First assume the special case that $M_1,\,M_2\in\matSU{n}{}$.
Then $X\in\ker L$ means $[X,\,M_1]=0$ and $[X,\,M_2]=0$.
Since $M_1$ and $M_2$ are irreducible,
by Lemma~\ref{lem:schur},
\begin{longversion}(lemma~\ref{lem:isid},)\end{longversion}
$X=r\id$ for some $r\in\bbC^\ast$.
Since $X$ is Hermitian, then $r\in\bbR^\ast$.
Hence $X$ is positive or negative definite.

To show the general case, let $C\in\matSL{n}{\bbC}$
be a simultaneous unitarizer of $M_1$ and $M_2$,
and let $\widetilde{L}=\calL(CM_1C^{-1},\,CM_2C^{-1})$.
Then $\widetilde{X} = {(C^{-1})}^\ast X C^{-1}$
is Hermitian,
and $\widetilde{X}\in\ker \widetilde{L}$
by Lemma~\ref{lem:L}(iii).
By the special case above, $\widetilde{X}$ is positive or negative definite.
Hence $X$ is positive or negative definite.
\end{proof}

\begin{lemma}
\label{lem:kernel}
Let $\calC$ be an analytic curve,
and $M_1, \dots ,M_q\in\Holo{\calC}{\matSL{n}{\bbC}}$ ($q\ge 2$),
and let $L = \calL(M_1,\dots,M_q)$.
Suppose for some subset $S\subset\calC$ of isolated points,
$M_1$ and $M_2$ are irreducible on $\calC\setminus S$,
and $M_1,\dots,M_q$ are pointwise simultaneously unitarizable
on $\calC\setminus S$.
Then there exists
$X\in(\Holo{\calC}{\matM{n\times n}{\bbC}}) \cap \ker L$ with $X^\ast=X$
and a subset $S'\subset\calC$ of isolated points such that
$X$ is positive definite on $\calC\setminus S'$.
\end{lemma}

\begin{proof}
By Lemma~\ref{lem:dimkerL},
$\dim\ker L=1$ on $\calC\setminus S$.
By Lemma~\ref{lem:holo-kernel},
$\dim\ker L\ge 1$ on $\calC$ and there exists
an analytic map $X_1\in\Holo{\calC}{\matM{n\times n}{\bbC}}$
such that $X_1\in(\ker L)\setminus\{0\}$.

If $X_1$ is Hermitian, define $X_2=X_1$;
otherwise define $X_2=i(X_1-X_1^\ast)$.
Then $X_2\not\equiv 0$, $X_2$ is Hermitian, and
$X_2\in\ker L$ by Lemma~\ref{lem:L}(ii).

By Lemma~\ref{lem:posdef},
$X_2$ is (pointwise) either positive definite or negative definite
except at the set of isolated points
at which $X_2=0$ or where $M_1,\dots,M_q$ all commute.

Let $v\in\bbC^n\setminus\{0\}$ be any non-zero constant vector
and define $f = v^\ast X_2 v$.
Then $f\not\equiv 0$ and
$X=f X_2$ is positive definite
on $\calC$ away from a set of isolated points.
\end{proof}

\typeout{===cholesky==============}
\Subsection{Analytic Cholesky decompositions}

We prove two analytic versions of the Cholesky decomposition theorem.
The first (Proposition~\ref{lem:cholesky})
is for Hermitian positive definite maps on an analytic curve,
used in the Unitarization Theorem~\ref{thm:1uni}.
The second (Proposition~\ref{lem:cholesky-semi}),
used in the $r$-unitarization Theorem~\ref{thm:runi},
is for meromorphic maps on $\bbS^1$
which are Hermitian positive definite except
at a finite set of points.

\begin{longversion}
First a preliminary lemma.

\begin{lemma}
\label{lem:cholesky-det}
Let $X\in\matM{n\times n}{\bbC}$ and write
\[
X = \begin{pmatrix} A & B \\ C & d\end{pmatrix},
\]
where
$A\in\matM{(n-1)\times(n-1)}{\bbC}$, $B\in\matM{(n-1)\times 1}{\bbC}$,
$C\in\matM{1\times(n-1)}{\bbC}$, and $d\in\bbC$.
If $A$ is invertible, then
\[
\frac{\det X}{\det A} = d - C A^{-1} B.
\]
\end{lemma}

\begin{proof}
\[
\begin{pmatrix}\id & 0 \\ -C & 1\end{pmatrix}
\begin{pmatrix}A^{-1} & 0 \\ 0 & 1\end{pmatrix}
\begin{pmatrix}A & B \\ C & d\end{pmatrix}
=
\begin{pmatrix}\id & A^{-1}B \\ 0  & d - CA^{-1}B\end{pmatrix}
\]
Taking the determinant yields the result.
\end{proof}
\end{longversion}

\begin{proposition}[Holomorphic Cholesky decomposition]
\label{lem:cholesky}
Let $\calC$ be an analytic curve
and let $X\in\Holo{\calC}{\matSL{n}{\bbC}}$ be
Hermitian positive definite.
Then
\begin{enumerate}
\item
There exists 
$V\in\Holo{\calC}{\matSL{n}{\bbC}}$ such that $X=V^\ast V$.
\item
$V$ is unique up to left multiplication by elements
of $\Holo{\calC}{\matSU{n}{}}$.
\end{enumerate}
\end{proposition}

\begin{proof}
We first prove the following
analytic version of the LDU-decomposition:
if $X\in\Holo{\calC}{\matGL{n}{\bbC}}$ is Hermitian
positive definite, then there exists
$R,\,D\in\Holo{\calC}{\matGL{n}{\bbC}}$
such that $X=R^\ast D R$,
where $R$ is upper triangular with diagonal elements $\equiv 1$,
and $D$ is diagonal
with diagonal elements in $\Holo{\calC}{\bbR^{>0}}$.
The proof is by induction on $n$.
The case $n=1$ is clear, with $R=\id$, $D=X$.

Now assume the statement is true for $n-1$ and write
\[
X = \begin{pmatrix} X_0 & Y_0 \\ Y_0^\ast & z\end{pmatrix},
\]
with $X_0\in\Holo{\calC}{\matGL{n-1}{\bbC}}$,
$Y_0\in\Holo{\calC}{\matM{(n-1)\times 1}{\bbC}}$ and
$z\in\Holo{\calC}{\bbC}$.
Then $X_0$ is Hermitian positive definite, so
let $X_0=R_0^\ast D_0 R_0$ be the decomposition
given by the induction hypothesis.
Then $D_0$ and $R_0^\ast$ are invertible on $\calC$, so
we can define $R,\,D\in\Holo{\calC}{\matGL{n}{\bbC}}$ by
\[
R = \begin{pmatrix} R_0 & S_0 \\ 0 & 1\end{pmatrix},\quad
S_0 = {(R_0^\ast D_0)}^{-1}Y_0,\quad
D = \begin{pmatrix}
D_0 & 0 \\ 0 & d
\end{pmatrix},\quad
d=z - S_0^\ast D_0 S_0.
\]
Then we have $X = R^\ast D R$.
Taking the determinant yields $\det X = \det D = d \det D_0$,
showing that $d$ takes values in $\bbR^{>0}$.
This proves the statement for $n$.

To show the first part of the theorem,
take Hermitian positive definite $X\in\Holo{\calC}{\matSL{n}{\bbC}}$,
and let $X=R^\ast D R$ be its analytic LDU factorization.
With $D=\diag(d_1,\dots,d_n)$, let
$E=\diag(\sqrt{d_1},\dots,\sqrt{d_n})$,
choosing positive square roots.
Let $V=ER$, so $X = V^\ast V$.
It is clear that $\det V\equiv 1$,
so $V\in\Holo{\calC}{\matSL{n}{\bbC}}$, proving (i).

To show uniqueness (ii), suppose
$V,\,W\in\Holo{\calC}{\matSL{n}{\bbC}}$ with
$V^\ast V = W^\ast W$.
Let $U=WV^{-1}\in\Holo{\calC}{\matSL{n}{\bbC}}$.
Then $U^\ast = U^{-1}$, so $U\in\Holo{\calC}{\matSU{n}{}}$.
\end{proof}

\begin{lemma}
\label{lem:unique-hermitian}
If $X_1,\,X_2\in\matSL{n}{\bbC}$ are positive definite,
and $X_2$ is a multiple of $X_1$,
then $X_1= X_2$.
\end{lemma}

\begin{proof}
Let $X_2 = c X_1$.
Since $X_1$ and $X_2$ are positive definite, then $c>0$.
Taking the determinant, $c^n=1$,
so $c$ is an $n$'th root of $1$.
Hence $c=1$, so $X_1 = X_2$.
\end{proof}

\begin{longversion}
\begin{lemma}
\label{lem:X-extender}
Let $\calC = (-1,\,1)$ and $X\in\Holo{\calC}{\matM{n\times n}{\bbC}}$.
Then
\begin{enumerate}
\item
If $X$ is Hermitian on $\calC\setminus\{0\}$, then
$X$ is Hermitian at $p$.
\item
If $X$ is positive definite on $\calC\setminus\{0\}$,
then $X$ is positive definite or positive semidefinite on $\calC$.
\item
If $X$ is positive definite on $(-1,\,0)$
and negative definite on $(0,\,1)$, then
$X(0)=0$.
\end{enumerate}
\end{lemma}
\end{longversion}

\typeout{===uni==============}
\Subsection{Simultaneous unitarization}

We are now prepared to prove the following unitarization theorem.
Conditions equivalent to condition (ii) of this theorem
(local simultaneous unitarizability) are found in
Sections~\ref{sec:graph}--\ref{sec:twobytwo}.

\begin{theorem}
\label{thm:1uni}
Let $\calC$ be an analytic curve and
$M_1,\dots,M_q\in\Holo{\calC}{\matSL{n}{\bbC}}$ ($q\ge 2$).
Suppose $M_1$ and $M_2$ are irreducible on $\calC$ except
at a subset of isolated points.
Then the following are equivalent:
\begin{enumerate}
\item
$M_1,\dots,M_q$ are globally simultaneously unitarizable on $\calC$.
\item
$M_1,\dots,M_q$ are locally simultaneously unitarizable
at each $p\in\calC$.
\end{enumerate}
In this case, any simultaneous unitarizer $V$
is unique up to left multiplication
by an element of $\Holo{\calC}{\matSU{n}{}}$.
\end{theorem}

\begin{proof}
Clearly (i) implies (ii).
Conversely, suppose (ii).
A global analytic simultaneous unitarizer $V$ is constructed as follows.


Fix $p\in\calC$. By (ii), there exists a neighborhood $\calU_p$
of $p$ and $W\in\Holo{\calU}{\matSL{n}{\bbC}}$ such that
$WM_jW^{-1}\in\Holo{\calU_p}{\matSU{n}{}}$, $j\in\{1,\dots,q\}$.
Define $X_p = W^\ast W$.
Then $X_p\in\Holo{\calU_p}{\matSL{n}{\bbC}}\cap\ker L$
and $X_p$ is Hermitian positive definite.
By Lemma~\ref{lem:unique-hermitian},
$X_p$ is the unique map in a neighborhood of $p$
which takes values in $\matSL{n}{\bbC}$,
is in $\ker L$,
and is Hermitian positive definite.

Define $X:\calC\to\matSL{n}{\bbC}$ by $X(p) = X_p(p)$.
Then for each $p\in\calC$, $X$ is analytic at $p$
because it coincides with the local analytic map $X_p$
on $\calU_p$, again by Lemma~\ref{lem:unique-hermitian}.
This defines the unique map
$X\in\Holo{\calC}{\matSL{n}{\bbC}}\cap\ker L$
which is Hermitian positive definite on $\calC$.

By the Cholesky Decomposition Proposition~\ref{lem:cholesky},
there exists
$V\in\Holo{\calC}{\matSL{n}{\bbC}}$ such that $X=V^\ast V$.
Then for $k\in\{1,\dots,q\}$,
$VM_kV^{-1}\in\Holo{\calC}{\matSU{n}{}}$
by Lemma~\ref{lem:L}(i),
so $V$ is the required simultaneous unitarizer.

To show uniqueness,
suppose $W\in\Holo{\calC}{\matSL{n}{\bbC}}$ is another
simultaneous unitarizer. Then $W^\ast W$
is a Hermitian positive definite
element of $\Holo{\calC}{\matSL{n}{\bbC}}\cap\ker L$.
By Lemma~\ref{lem:unique-hermitian}, $W^\ast W=X=V^\ast V$.
The uniqueness result follows by the uniqueness result in
the Cholesky Decomposition Proposition~\ref{lem:cholesky}.
\end{proof}

\Section{Local conditions for simultaneous unitarizability}

Theorem~\ref{thm:1uni} effectively reduces
global simultaneous unitarizability
to local simultaneous unitarizability.
We will now give more explicit conditions for
local simultaneous unitarizability.

\typeout{===graph==============}
\Subsection{$\Xi_n$-graphs}
\label{sec:graph}

Given two matrix maps $M_1,\,M_2\in\Holo{\calC}{\matSL{n}{\bbC}}$
with $M_1$ diagonal,
the local simultaneous unitarizability of $M_1$ and $M_2$ at $p$
is equivalent to the equalities of the orders of
certain corresponding entries of
$M_2$ and ${M_2^\ast}^{-1}$ at $p$.
These relations are naturally expressed
in terms of graphs on $\{1,\dots,n\}$ (Definition~\ref{def:Qgraph}).


\begin{definition}
\label{def:Qgraph}
A \emph{$\Xi_n$-graph}
is a directed graph with vertices $V = \{1,\dots,n\}$.
A $\Xi_n$-graph is not a multigraph
--- it may not have two instances
of the same edge.
However, it may have an edge connecting a vertex
to itself. A $\Xi_n$-graph is \emph{connected}
if it is connected as a undirected graph.

Let $\calC$ be an analytic curve,
let $p\in\calC$,
and let $A,\,B\in\Holo{\calC}{\matM{n\times n}{\bbC}}$.
For $\mu,\nu\in\{1,\dots,n\}$, let $A_{\mu\nu}$ denote the entry
of $A$ lying in row $\mu$ and column $\nu$, and similarly for $B$.

Let $G$ be a $\Xi_n$-graph.
$A$ is \emph{$G$-non-zero at $p$}
if for every directed edge $(\mu,\,\nu)\in V^2$
from $\mu$ to $\nu$ of $G$
we have $\ord_p A_{\mu\nu}<\infty$.
(Since $\calC$ is connected and $A$ is holomorphic,
this order condition is equivalent to the condition $A_{\mu\nu}\not\equiv 0$.)

Let $G$ be a $\Xi_n$-graph.
$A$ and $B$ are
\emph{$G$-compatible at $p$}
if for every directed edge $(\mu,\,\nu)\in V^2$ of $G$,
we have $\ord_p A_{\mu\nu} = \ord_p B_{\mu\nu} < \infty$.
\end{definition}

Let $e_1=(1,\,0,\dots,0),\dots,e_n=(0,\dots,0,\,1)\in\bbC^n$
be the standard basis for $\bbC^n$.
%
Fixing $n$, we have the natural correspondence between
the nonempty subsets of $\{1,\dots,n\}$ and the
nonzero subspaces of $\bbC^n$ generated by standard basis elements,
where $K\subset\{1,\dots,n\}$ corresponds to the subspace
\[
E_K = \Span\{e_k\suchthat k\in K\}\subseteq\bbC^n.
\]
This correspondence induces a
natural one-to-one correspondence between
the set of partitions of $\{1,\dots,n\}$
all of whose terms are non-empty, and the set of
direct sum decompositions of $\bbC^n$ relative to its standard basis,
all of whose terms are non-zero.

The following lemma expresses
in terms of $\Xi_n$-graphs
a notion closely related to block-diagonalizability.

\begin{lemma}
\label{lem:Qgraph1}
Let $\calC$ be an analytic curve,
let $A\in\Holo{\calC}{\matM{n\times n}{\bbC}}$,
and let $p\in\calC$.
Then the following are equivalent:
\begin{enumerate}
\item
Every
$\Xi_n$-graph $H$ for which $A$ is $H$-non-zero at $p$
is disconnected.
\item
There exists a direct decomposition
$\bbC^n = V_1\oplus V_2$ into non-zero summands
$V_1$ and $V_2$
generated by standard basis elements of $\bbC^n$
such that $AV_1\subseteq V_1$ and $AV_2\subseteq V_2$ near $p$.
\end{enumerate}
\end{lemma}

\begin{proof}
Suppose (ii) and let $G_1\sqcup G_2 = \{1,\dots,n\}$ be the corresponding
partition.
The conditions $AV_1\subseteq V_1$ and $AV_2\subseteq V_2$
are equivalent to the condition that
$A_{\mu\nu}\equiv 0$ for all
$(\mu,\,\nu)\in(G_1\times G_2)\cup(G_2\times G_1)$.
Let $G$ be a $\Xi_n$-graph for which $A$ is $G$-non-zero at $p$.
Then $G$ does not contain any of the edges
in $(G_1\times G_2)\cup(G_2\times G_1)$.
Hence the two subsets $G_1$ and $G_2$ of the vertex set of $G$
are disconnected, proving (i).

Conversely, suppose (i) and
let $G$ be the maximal $\Xi_n$-graph among the $\Xi_n$-graphs
$H$ for which $A$ is $H$-non-zero at $p$.
Then $G$ is disconnected, so
let $G_1\sqcup G_2=\{1,\dots,n\}$ be a partition
with $G_1$ and $G_2$ nonempty
such that no edge of $G$ is in
$(G_1\cross G_2)\cup(G_2\cross G_1)$.
Since $G$ is maximal, then $A_{\mu\nu}\equiv 0$
for all $(\mu,\,\nu)\in (G_1\times G_2)\cup(G_2\times G_1)$.
Let $\bbC^n = V_1\oplus V_2$ be the direct sum decomposition
corresponding to the partition $G_1\sqcup G_2$.
Then $AV_1\subseteq V_1$ and $AV_2\subseteq V_2$,
proving (ii).
\end{proof}

We shall say that $X\in\Holo{\calC}{\matM{n\times n}{\bbC}}$ is
\emph{infinitesimally invertible at $p$} if
the leading term in its series expansion at $p$
in some local coordinate on $\calC$ near $p$ has full rank.
Equivalently,
there exists a local meromorphic scalar function $f$ near $p$
such that $fX$ is holomorphic at $p$ and $\rank (fX)(p)=n$.

\begin{lemma}
\label{lem:Qgraph2}
Let $\calC$ be an analytic curve,
let $p\in\calC$
and let $A,\,B\in\Holo{\calC}{\matM{n\times n}{\bbC}}$.
Let $X\in\Holo{\calC}{\matM{n\times n}{\bbC}}$ be diagonal
with $X\not\equiv 0$,
and suppose
$X A = B X$.
Then the following are equivalent:
\begin{enumerate}
\item
$A$ and $B$ are $G$-compatible at $p$ for some
connected $\Xi_n$-graph $G$.
\item
$X$ is infinitesimally invertible at $p$,
and neither of the equivalent conditions of Lemma~\ref{lem:Qgraph1}
hold for $A$ at $p$.
\end{enumerate}
\end{lemma}

\begin{proof}
First suppose (ii) holds.
Since $X$ is diagonal,
the infinitesimal invertibility of $X$ is equivalent to
\begin{equation}
\label{eq:Xconj-1}
\ord_p X_{11} = \dots = \ord_p X_{nn} < \infty.
\end{equation}
The components of $XA-BX$ are
\begin{equation}
\label{eq:Xconj-2}
0 = (XA-BX)_{\mu\nu} = X_{\mu\mu} A_{\mu\nu} - X_{\nu\nu} B_{\mu\nu}.
\end{equation}
Let $G$ be a connected $\Xi_n$-graph such that
$A$ is $G$-non-zero at $p$.
Then for each edge $(\mu,\,\nu)$ of $G$,
we have $\ord_p A_{\mu\nu} < \infty$.
Then~\eqref{eq:Xconj-1} and~\eqref{eq:Xconj-2} imply
\begin{equation}
\label{eq:Xconj-3}
\ord_p A_{\mu\nu} = \ord_p B_{\mu\nu} < \infty.
\end{equation}
Hence $A$ and $B$ are $G$-compatible at $p$.

Conversely, suppose (i), that $A$ and $B$ are $G$-compatible at $p$
for some connected $\Xi_n$-graph $G$.
Then~\eqref{eq:Xconj-3} holds for each edge $(\mu,\,\nu)$ of $G$.
Hence $A$ is $G$-non-zero at $p$.
By~\eqref{eq:Xconj-2},
\[
\ord_p X_{\mu\mu} = \ord_p X_{\nu\nu} < \infty.
\]
The connectedness of $G$ implies~\eqref{eq:Xconj-1}.
\end{proof}

\begin{lemma}
\label{lem:detX}
Let $\calC$ be an analytic curve,
$p\in\calC$,
and suppose $X\in\Holo{\calC}{\matM{n\times n}{\bbC}}$
is infinitesimally invertible at $p$.
Then
\begin{enumerate}
\item
$\det X$ has a local analytic $n$'th root in some neighborhood
$\calU\subset\calC$ of $p$.
\item
${(\det X)}^{-1/n}X\in\Holo{\calU}{\matSL{n}{\bbC}}$.
\item
If $X$ is Hermitian positive definite on a punctured
neighborhood of $p$, then the $n$'th root can
be chosen so that
${(\det X)}^{-1/n}X$ is Hermitian positive definite
in a neighborhood of $p$.
\end{enumerate}
\end{lemma}

\begin{proof}
Let $t$ be a local coordinate on $\calC$ near $p$ with
$t(p)=0$. 
Since $X$ is infinitesimally invertible at $p$, we have
\begin{equation}
\label{eq:detX1}
X = X_m t^m + \Order(t^{m+1}),\quad \det X_m\ne 0.
\end{equation}
Taking the determinant of~\eqref{eq:detX1},
\begin{equation}
\label{eq:detX2}
\det X = (\det X_m) t^{nm} + \Order(t^{nm+1}).
\end{equation}
Since $\det X_m\ne 0$, then $\ord_p(\det X) = nm$.
Since $nm$ is divisible by $n$, then
$\det X$ has a local analytic $n$'th root
in a neighborhood $\calU\subset\calC$ of $p$, proving (i).

Taking an $n$'th root of~\eqref{eq:detX2},
\begin{equation}
\label{eq:detX3}
{(\det X)}^{1/n} = {(\det X_m)}^{1/n} t^{m} + \Order(t^{m+1}),
\end{equation}
so
\[
Y:={(\det X)}^{-1/n} X = {(\det X_m)}^{-1/n} X_m + \Order(t^1),
\]
and hence
$Y\in\Holo{\calU}{\matSL{n}{\bbC}}$,
proving (ii).

To show (iii),
suppose $\calU$ is a neighborhood of $p$ such that
$X$ is positive definite on $\calU\setminus\{p\}$.
Then $\det X$ is positive on $\calU\setminus\{p\}$,
and hence by its continuity is nonnegative on $\calU$.
By~\eqref{eq:detX2}, $\det X_m>0$, so
we can choose ${(\det X_m)}^{1/n}>0$ in~\eqref{eq:detX3}.
Equation~\eqref{eq:detX1} and a limit argument imply that $m$ is even.
Equation~\eqref{eq:detX3} then implies that ${(\det X)}^{1/n}$ is non-negative
on $\calU$.
Then $Y={(\det X)}^{1/n}X$
is positive definite in a punctured neighborhood of $p$.
Hence $Y(p)$ is positive semidefinite, and so $Y(p)$ is positive
definite since $\det Y(p)=1$.
\end{proof}

\begin{lemma}
\label{lem:uni-Qgraph}
Let $\calC$ be an analytic curve,
let $p\in\calC$ and
let $M_1,\,M_2\in\Holo{\calC}{\matSL{n}{\bbC}}$.
Suppose $M_1$ and $M_2$ are irreducible on $\calC\setminus\{p\}$.
Suppose $M_1$ is diagonal
and no two local analytic eigenvalues of $M_1$ are identically equal.
Then the following are equivalent:
\begin{enumerate}
\item
$M_1$ and $M_2$ are locally simultaneously unitarizable at $p$.
\item
$M_1$ and $M_2$ are pointwise simultaneously unitarizable at
each point of a neighborhood of $p$, and
$M_2$ and ${M_2^\ast}^{-1}$ are $G$-compatible at $p$ for
some connected $\Xi_n$-graph $G$.
\end{enumerate}
\end{lemma}

\begin{proof}
To show (i) $\Rightarrow$ (ii),
let $V\in\Holo{\calU}{\matSL{n}{\bbC}}$
be a local simultaneous unitarizer of $M_1$ and $M_2$ at $p$
and let $X = V^\ast V$.
Then
\[
XM_k = {M_k^\ast}^{-1}X,\quad k\in\{1,\,2\}.
\]
Since $M_1$ is unitarizable, it has unimodular
eigenvalues. Since $M_1$ is diagonal, then $M_1\in\Holo{\calC}{\matSU{n}{}}$.
Hence $M_1 = {M_1^\ast}^{-1}$, and so $[X,\,M_1]=0$.
Since no two eigenvalues of $M_1$ are identically equal,
$X$ is diagonal away from the set of isolated points where two
eigenvalues of $M_1$ coincide, so by its continuity,
$X$ is diagonal.

Since $M_1$ is diagonal and
$M_1$ and $M_2$ are irreducible in a punctured neighborhood of $p$, then
there is no non-zero proper
subspace $W\subset\bbC^n$ generated by standard
basis elements of $\bbC^n$ such that $M_2W\subset W$
in a neighborhood of $p$.
Hence $M_2$ is not block-diagonalizable at $p$ via a permutation at $p$.
Since $X(p)\in\matSL{n}{\bbC}$, then $X$ has full rank at $p$.
By Lemma~\ref{lem:Qgraph2},
$M_2$ and ${M_2^\ast}^{-1}$ are $G$-compatible at $p$
for some connected $\Xi_n$-graph $G$.
This proves (ii).

To show (ii) $\Rightarrow$ (i),
since $M_1$ and $M_2$ are irreducible on $\calC\setminus\{p\}$
and are pointwise simultaneous unitarizable at each point in
a neighborhood of $p$,
by Lemma~\ref{lem:kernel}
there exists a neighborhood $\calU$ of $p$ and a map
$X\in\Holo{\calU}{\matM{n\times n}{\bbC}}\cap\ker L$
for which $X$ is Hermitian positive definite on $\calU\setminus\{p\}$.
By condition (ii) and Lemma~\ref{lem:Qgraph2},
$X$ is infinitesimally invertible at $p$.
By Lemma~\ref{lem:detX},
there exists a neighborhood $\calV$ of $p$
and a choice of $n$'th root of $\det X$ such that
$Y = {(\det X)}^{-1/n}X\in\Holo{\calV}{\matSL{n}{\bbC}}$
is Hermitian positive definite.

By the Cholesky Decomposition Proposition~\ref{lem:cholesky},
there exists
$V\in\Holo{\calV}{\matSL{n}{\bbC}}$ such that $Y=V^\ast V$.
Then by Lemma~\ref{lem:L}(i),
$VM_kV^{-1}\in\Holo{\calV}{\matSU{n}{}}$, $k\in\{1,\,2\}$,
so $V$ is a local simultaneous unitarizer of $M_1$ and $M_2$ at $p$.
\end{proof}

\begin{lemma}
\label{lem:free-uni}
For $q\ge 2$, let $M_1,\dots,M_q\in\matSL{n}{\bbC}$
with $M_1,\,M_2$ irreducible
and $M_1,\,M_2\in\matSU{n}{}$.
If $M_1,\dots,M_q$ are simultaneously unitarizable
by an element of $\matSL{n}{\bbC}$,
then $M_1,\dots,M_n\in\matSU{n}{}$.
\end{lemma}

\begin{proof}
Let $C\in\matSL{n}{\bbC}$ be a simultaneous unitarizer.
Then $C^\ast C$ commutes with each of $M_1$ and $M_2$.
By Lemma~\ref{lem:schur},
$C^\ast C\in\bbC\id$.
Since $C^\ast C$ is positive definite,
by Lemma~\ref{lem:unique-hermitian}, $C^\ast C=\id$.
Hence $C\in\matSU{n}{}$.
Since $C M_k C^{-1}\in\matSU{n}{}$ for $k\in\{1,\dots,q\}$,
then $M_k\in\matSU{n}{}$ for $k\in\{1,\dots,q\}$.
\end{proof}

Lemmas~\ref{lem:uni-Qgraph}, \ref{lem:free-uni}
and Theorem~\ref{thm:1uni} give the following
result.

\begin{theorem}[Unitarization theorem]
\label{thm:nuni}
Let $\calC$ be an analytic curve,
let $p\in\calC$ and
let $M_1,\dots,M_q\in\Holo{\calC}{\matSL{n}{\bbC}}$, $q\ge 2$.
Suppose $M_1$ and $M_2$ are irreducible on $\calC$ except
at a subset of isolated points.
Suppose $M_1$ is locally diagonalizable at each point $p\in\calC$,
and that 
no two local analytic eigenvalues of $M_1$ are identically equal.
Then
$M_1,\dots,M_q$ are globally simultaneously unitarizable
if and only if the following conditions hold:
\begin{enumerate}
\item
$M_1,\dots,M_q$ are pointwise simultaneously unitarizable at
each point of $\calC$.
\item
For each $p\in\calC$, let $C$ be a local diagonalizer of $M_1$
at $p$, and let $P_2 = CM_2 C^{-1}$.
Then $P_2$ and ${P_2^\ast}^{-1}$ are $G$-compatible at $p$ for
some connected $\Xi_n$-graph $G$.
\end{enumerate}
In this case, any simultaneous unitarizer is unique
up to left multiplication
by an element of $\Holo{\calC}{\matSU{n}{}}$.
\end{theorem}

\Section{Simultaneous unitarizability for $\matSL{2}{\bbC}$}
\label{sec:twobytwo}

\typeout{===twobytwo==============}

For the case $\matSL{2}{\bbC}$, the $\Xi_n$-graph condition
in Theorem~\ref{thm:1uni}
is particularly simple and can be recast in terms of
commutators.

Note that for $A,\,B\in\matM{2\times 2}{\bbC}$,
$A$ and $B$ are irreducible if and only if $[A,\,B]$ has rank 2.

\begin{definition}
\label{def:infinitesimally-irreducible}
We say that $A,\,B\in\Holo{\calC}{\matSL{2}{\bbC}}$ are
\emph{infinitesimally irreducible} at $p\in\calC$
if $[A,\,B]\not\equiv 0$ and the leading
order term in the series expansion of $[A,\,B]$ at $p$ has full rank.
\end{definition}

The property infinitesimally irreducible
(resp.~reducible) is preserved
under conjugation by an element of $\Holo{\calC}{\matSL{2}{\bbC}}$.

\begin{lemma}
\label{lem:diagonal-infinitesimally-irreducible}
If $A$ is diagonal, then $A$ and $B$ are infinitesimally
irreducible if and only if the two off-diagonal terms of $B$
have the same finite order.
\end{lemma}

We give a sufficient condition for local diagonalizability:

\begin{lemma}
\label{lem:SUn-local-diagonalizability}
Let $\calC$ be an analytic curve, $p\in\calC$, and
$M\in\Holo{\calC}{\matSU{2}{}}$ be an analytic map.
Then $M$ is locally diagonalizable in a neighborhood $\calU\subset\calC$
of $p$ by a map $C\in\Holo{\calU}{\matSU{2}{}}$.
\end{lemma}

\begin{proof}
It follows from the characteristic equation of $M$ that
its eigenvalues $\mu_1,\,\mu_2 = \mu_1^{-1}$
are analytic in a neighborhood $\calU\subset\calC$ of $p$.
It can be shown, for example by an analytic version of the QR-decomposition,
 that there exist corresponding analytic eigenvector functions
$v_1,\,v_2\in\Holo{\calU}{\bbC^2}$ such that
$V = (v_1,\,v_2)\in\Holo{\calU}{\matSU{2}{}}$.
Then $V^{-1}MV = \diag(\mu_1,\,\mu_2)$, so $C=V^{-1}$ is the
required diagonalizer.
\end{proof}


\begin{theorem}[Unitarization theorem for $\matSL{2}{\bbC}$]
\label{thm:uni-2by2}
Let $\calC$ be an analytic curve,
$M_1, \dots ,M_q\in\Holo{\calC}{\matSL{2}{\bbC}}$ ($q\ge 2$),
and suppose that $[M_r,\,M_s]\not\equiv 0$
for some fixed choice $r,\,s\in\{1,\dots,q\}$.
Then $M_1,\dots,M_q$ are globally simultaneously unitarizable
if and only if the following conditions hold
at each $p\in\calC$:
\begin{enumerate}
\item
$M_1,\dots,M_q$ are pointwise simultaneously unitarizable at $p$.
\item
$M_r$ or $M_s$ is locally diagonalizable at $p$.
\item
$M_r$ and $M_s$ are infinitesimally irreducible
at $p$.
\end{enumerate}
In this case, any simultaneous unitarizer $V$ is unique
up to left multiplication
by an element of $\Holo{\calC}{\matSU{2}{}}$.
\end{theorem}

\begin{remark}
Examples exist which show the independence of the three conditions (i)--(iii)
of Theorem~\ref{thm:uni-2by2}.
\end{remark}

\begin{proof}
Renumber so that $r=1$ and $s=2$.
Suppose $M_1,\dots,M_q$ are globally simultaneously unitarizable on $\calC$.
Then (i) clearly holds, and
condition (ii) holds by Lemma~\ref{lem:SUn-local-diagonalizability}.
To show condition (iii), let $V$ be a local unitarizer of
$M_1$ and $M_2$, and let $P_k = V M_k V^{-1}$, $k\in\{1,\,2\}$.
By Lemma~\ref{lem:SUn-local-diagonalizability}, there exists a
local unitary diagonalizer $C$ of $P_1$.
Let $Q_k = C P_k C^{-1}$, $k\in\{1,\,2\}$.
Then $Q_1$ is diagonal.
Since $Q_2$ is unitary, its off-diagonal terms have the
same order.
Since $[M_1,\,M_2]\not\equiv 0$, then
$[Q_1,\,Q_2]\not\equiv 0$, so $Q_2$ is not identically diagonal.
Hence the off-diagonal entries of $Q_2$ are not both identically zero.
By Lemma~\ref{lem:diagonal-infinitesimally-irreducible},
$Q_1$ and $Q_2$ are
infinitesimally irreducible,
and hence $M_1$ and $M_2$ are infinitesimally irreducible.

Conversely, assume conditions (i)--(iii)
and assume $M_1$ is locally diagonalizable at $p$.
We show that the conditions of Theorem~\ref{thm:nuni} are
satisfied.
Since $M_1$ is locally diagonalizable, and $[M_1,\,M_2]\not\equiv 0$,
then $M_1$ does not have identically equal eigenvalues.
Let $C$ be a local diagonalizer of $M_1$ at $p$,
and let $P_k = C M_k C^{-1}$, $k\in\{1,\,2\}$.
By (iii), the two off-diagonal terms of $P_2$ and those of
${P_2^\ast}^{-1}$ all have the same finite order.
Hence $P_1$ and $P_2$ are irreducible in a punctured neighborhood of $p$,
so $M_1$ and $M_2$ are irreducible in a punctured neighborhood of $p$.
Let $G$ be the connected $\Xi_2$-graph with single edge $(1,\,2)$.
Then $P_2$ and ${{P_2}^\ast}^{-1}$ are $G$-compatible.
The existence and uniqueness
of the global simultaneous unitarizer follows by Theorem~\ref{thm:nuni}.
\end{proof}

\section{Simultaneous $r$-unitarization}
\label{sec:runi}

\typeout{===runi==============}

The $r$-Unitarization Theorem~\ref{thm:runi}
for $\matSL{n}{\bbC}$-valued loops
is a variant of the Unitarization Theorem~\ref{thm:1uni}
on the standard unit circle
$\bbS^1 = \{\lambda\in\bbC\suchthat\abs{\lambda}=1\}$.
This variant has been proven for the case of $\matSL{2}{\bbC}$
\cite{SKKR:spaceforms},
where it finds application to the construction of non-simply-connected
CMC surfaces.

In the $r$-Unitarization Theorem,
a holomorphic map $V$ is constructed on an annulus
$\calA_{r,1}$ which simultaneously unitarizes
the given set of loops $M_1,\dots,M_q$ in the sense
that the $VM_kV^{-1}$ extend
holomorphically to $\bbS^1$ and are unitary there.
In this case, the unitarizing loop $V$ does not in general extend
holomorphically to $\bbS^1$, but has
roots of zeros and poles there.

We note that a unitary map cannot have poles.
This follows from the fact that
$\matSU{n}{}$ and $\matU{n}{}$ are compact:

\begin{proposition}
\label{lem:uni-extender}
Let $\calC$ be an analytic curve. Then
$\Mero{\calC}{\matU{n}{}}=\Holo{\calC}{\matU{n}{}}$
and $\Mero{\calC}{\matSU{n}{}}=\Holo{\calC}{\matSU{n}{}}$.
\end{proposition}

\begin{longversion}
\begin{proof}
$\matSU{n}{}$ is bounded in the Euclidean vector norm $\norm{\,\cdot\,}$.
On the other hand, if $M\in\Mero{\calC}{\matM{n\times n}{\bbC}}$
has a pole at $p\in\calC$, $\norm{M}$ is unbounded in any
neighborhood of $p$.
Hence $U\in\Mero{\calC}{\matSU{n}{}}$ cannot have a pole on $\calC$.
\end{proof}
\end{longversion}

We now extend the Cholesky decomposition theorem to the
case of an analytic map on $\bbS^1$ which is Hermitian positive definite
except at a finite subset.

\begin{proposition}[Meromorphic Cholesky decomposition]
\label{lem:cholesky-semi}
Let $\calC=\bbS^1\subset\bbC$
and let $X\in\Mero{\calC}{\matM{n\times n}{\bbC}}$ be
Hermitian positive definite except at a finite subset of points
$S\subset\calC$. Then
\begin{enumerate}
\item
There exists 
$V\in\Mero{\calC}{\matM{n\times n}{\bbC}}$
such that $X=V^\ast V$.
\item
$V$ is unique up to left multiplication by elements
of $\Holo{\calC}{\matU{n}{}}$.
\end{enumerate}
\end{proposition}

\begin{proof}
We will apply the LDU-decomposition stated at the beginning of the proof of
Proposition~\ref{lem:cholesky},
with holomorphicity replaced by meromorphicity.

Let $\rho:\widetilde{\calC}\to\calC$ be a double cover
and let $\tau:\widetilde{\calC}\to\widetilde{\calC}$
be the deck transformation
induced by a single counterclockwise traversal of $\calC$.
Let $\rho^\ast $ and $\tau^\ast$ denote the respective pullbacks.
Write $D=\diag(d_1,\dots,d_n)$ and for $k\in\{1,\dots,n\}$,
define $b_k$ to be either of the global square roots of $d_k$
on $\widetilde\calC$.
Define $\widetilde{B}=\diag(b_1,\dots,b_n)$ and
$\widetilde{V} = \widetilde{B}(\rho^\ast R)$,
so $\rho^\ast X = \widetilde{V}^\ast \widetilde{V}$.
With $\lambda$ the standard unimodular parameter on $\bbS^1$,
define
\[
\text{$c_k = 1$ if $\tau^\ast b_k = b_k$,}
\quad\text{and}\quad
\text{$c_k = \sqrt{\lambda}$ if $\tau^\ast b_k = -b_k$,}
\]
and define
$C=\diag(c_1,\dots,c_n)\in\Holo{\widetilde{\calC}}{\matU{n}{}}$.
Then $\tau^\ast (C\widetilde{V}) = C\widetilde{V}$.
Let $V\in\Mero{\calC}{\matM{n\times n}{\bbC}}$
be the unique map satisfying $\rho^\ast V = U\widetilde{V}$.
Then $X = V^\ast V$, proving (i).

To show uniqueness (ii), suppose
$V,\,W\in\Mero{\calC}{\matM{n\times n}{\bbC}}$ with
$V^\ast V = W^\ast W$.
Let $U=WV^{-1}\in\Mero{\calC}{\matM{n\times n}{\bbC}}$.
Then $U^\ast = U^{-1}$, so $U$ takes values in $\matU{n}{}$ on
$\bbS^1$ away from its poles. By Proposition~\ref{lem:uni-extender},
$U\in\Holo{\calC}{\matU{n}{}}$.
\end{proof}

\begin{definition}
\label{def:r-unitarizable}
Let $M_1,\dots,M_q\in\Holo{\bbS^1}{\matSL{n}{\bbC}}$.
Let $r\in(0,\,1)$ and suppose
$M_1,\dots,M_q$ extend holomorphically to respective maps
$\widetilde{M}_1,\dots\widetilde{M}_q\in\Holo{\calA_{r,1}}{\matSL{n}{\bbC}}$.
Then $M_1,\dots,M_q$ are \emph{simultaneously $r$-unitarizable}
if there exists
$V\in\Holo{\calA_{r,1}}{\matSL{n}{\bbC}}$ which
extends holomorphically to $\bbS^1$ minus a finite subset such that
$V\widetilde{M}_1V^{-1},\dots,V\widetilde{M}_qV^{-1}$
extend holomorphically to $\bbS^1$
and their respective restrictions to $\bbS^1$
are in $\Holo{\bbS^1}{\matSU{n}{}}$.
\end{definition}

\begin{theorem}[$r$-unitarization theorem]
\label{thm:runi}
Let $M_1, \dots ,M_q\in\Holo{\bbS^1}{\matSL{n}{\bbC}}$ ($q\ge 2$)
and suppose for some $i,\,j\in\{1,\dots,r\}$ that
$M_i$ and $M_j$ are irreducible except at a finite subset of $\bbS^1$.
Then the following are equivalent:
\begin{enumerate}
\item
$M_1,\dots,M_q$ are simultaneously $r$-unitarizable
for some $r\in(0,\,1)$.
\item
$M_1,\dots,M_q$ are pointwise simultaneously unitarizable
on $\bbS^1$ minus a finite subset.
\end{enumerate}
Moreover, $r$-unitarizers are unique in the following sense.
If $V_1\in\Holo{\calA_{r_1,1}}{\matSL{n}{\bbC}}$
and $V_2\in\Holo{\calA_{r_2,1}}{\matSL{n}{\bbC}}$
are respective $r_1$- and $r_2$- unitarizers,
then $V_2V_1^{-1}$
extends holomorphically to $\bbS^1$ and
its restriction to $\bbS^1$ is in
$\Holo{\bbS^1}{\matSU{n}{}}$.
\end{theorem}

\begin{proof}
First suppose a simultaneous $r$-unitarizer $V$ exists as in (i)
and let $S\subset\bbS^1$ be the finite singular set of the extension of
$V$ to $\bbS^1$.
Then by the definition of $r$-unitarizer, for all $p\in\bbS^1\setminus S$,
$\left.VM_kV^{-1}\right|_p\in\matSU{n}{}$, proving (ii).

Conversely, suppose (ii) holds.
A simultaneous $r$-unitarizer $V$ is constructed as follows.

Let $L$ be as in Definition~\ref{def:L}.
By Lemma~\ref{lem:kernel},
there exists $X_1\in\Holo{\calC}{\matM{n\times n}{\bbC}}$ such that
$X_1\in\ker L$,
$X_1^\ast = X_1$,
and away from a finite subset of $\bbS^1$,
$X_1$ is positive definite.

By the Cholesky Decomposition Proposition~\ref{lem:cholesky-semi}
there exists $V_1\in\Mero{\calC}{\matM{n\times n}{\bbC}}$
such that $X_1=V_1^\ast V_1$.

By the holomorphicity of $M_1,\dots,M_q$ and the meromorphicity
of $V_1$, there exists $r\in(0,\,1)$
such that $M_1,\dots,M_q$ and $V_1$ extend holomorphically to $\calA_{r,1}$
and $\det V_1$ is non-zero in $\calA_{r,1}$.

Let $\widetilde{\calA}\to\calA_{r,1}$ be an $n$-fold cover
and let $\tau:\widetilde{\calA}\to\widetilde{\calA}$
be the deck transformation induced by
a single counterclockwise traversal of $\bbS^1$.
Define $V_2 = {(\det V)}^{-1/n}V_1$ on $\widetilde{\calA}$.
Then for some $k\in\bbZ^{\ge 0}$, $\tau^\ast V_2 = \epsilon^k V_2$,
where $\epsilon=e^{2\pi i/n}$.
Let $\lambda$ be the standard unimodular parameter on $\bbS^1$.
Define $U\in\Holo{\widetilde{\calA}}{\matSU{n}{}}$ by
\[
U = \lambda^{-k/n}\diag(1,\dots,1,\lambda^{k}),
\]
so $\tau^\ast U = \epsilon^{-k}U$.
Let $V = UV_2$.
Then $\tau^\ast V = V$,
so $V$ on $\widetilde{\calA}$ descends to a single valued
holomorphic map on $\calA_{r,1}$.
This gives us $V\in\Holo{\calA_{r,1}}{\matSL{2}{\bbC}}$.

Let $S\subset\bbS^1$ be the singular set of $V$.
By Lemma~\ref{lem:L}(i),
on $\bbS^1\setminus S$
for each $k\in\{1,\dots,q\}$,
$P_k:=V M_k V^{-1}$ takes values in $\matSU{n}{}$.
Since $V_1$ is meromorphic on $S$,
then $P_k = V_1 M_k V_1^{-1}$ is meromorphic on $S$.
By Proposition~\ref{lem:uni-extender},
$P_k\in\Holo{\bbS^1}{\matSU{n}{}}$.
Hence $V$ is the required simultaneous $r$-unitarizer.

The uniqueness result follows as in the proof of Theorem~\ref{thm:1uni}.
\end{proof}

\typeout{===weierstrass=============}

\Section{The extended Weierstrass representation}
\label{sec:weierstrass}

We now construct trinoids and symmetric $n$-noids,
conformal CMC immersions of the $n$-punctured Riemann sphere
into each of the space forms
Euclidean $3$-space $\bbR^3$, spherical $3$-space $\bbS^3$ and
hyperbolic $3$-space $\bbH^3$.
We first describe the extended
Weierstrass representation used in the construction.
As it is thoroughly described in a number of places,
such as~\cite{DPW,DorH:cyl,SKKR:spaceforms}, we give only a
brief outline here.

\Subsection{The Iwasawa decomposition}

Given an analytic Lie group $G$,
we denote by $\Lambda G$
the group $C^\omega(\bbS^1,\,G)$ of analytic maps $\bbS^1\to G$.
Let $\calD_1 \subset \bbC$ be the open disk bounded by $\bbS^1$.
The subgroup $\LooppSL{}\subset\LoopSL{}$ of \emph{positive loops}
is the subgroup of loops $B\in\LoopSL{}$
such that $B$ extends holomorphically to $\calD_1$ and $B(0)$ is
upper triangular with real diagonal entries.
The subgroup $\LoopuSL{}\subset\LoopSL{}$ of \emph{unitary loops}
is the subgroup of loops $F\in\LoopSL{}$ which satisfy the
condition
$F^\ast = F^{-1}$,
where for any $F\in\LoopSL{}$, $F^\ast\in\LoopSL{}$ is defined by
\begin{equation}
\label{eq:star}
F^\ast(\lambda) = \transpose{\overline{F(1/\ol{\lambda})}} \spaceperiod
\end{equation}
Note that $F\in \LoopuSL{}$ implies $F(p) \in \matSU{2}{}$
at each point $p\in\bbS^1$.

Multiplication $\LoopuSL{} \times \LooppSL{} \to \LoopSL{}$
is a real-analytic diffeomorphism onto~\cite{Pressley-Segal,DPW}.
For $\Phi \in \LoopSL{}$,
\begin{equation*}
  \Phi = FB \spacecomma
\end{equation*}
with $F \in \LoopuSL{}$ and $B \in \LooppSL{}$,
is called the $1$-Iwasawa (or just Iwasawa) decomposition of $\Phi$.
The chosen normalization of $B(0)$ gives
uniqueness of this decomposition.  We call $F$ the unitary factor of $\Phi$.

\Subsection{The extended Weierstrass construction}

Every conformal CMC $H$ immersion into one of the $3$-dimensional
space forms $\bbR^3$ or $\bbS^3$ or $\bbH^3$
can be locally constructed by
the extended Weierstrass representation~\cite{SKKR:spaceforms,DPW} (with 
$H \neq 0$ for $\bbR^3$ and $|H|>1$ for $\bbH^3$) 
as follows:

1. Let $\Sigma$ be a domain in the $z$-plane $\bbC$, and choose a
holomorphic
$C^\omega(\bbS^1,\,\matsl{2}{\bbC})$-valued
differential form $\xi = A(z,\lambda)dz$
which extends meromorphically to $\calD_1$ with a pole only at $\lambda=0$,
which is simple and appears only in the upper-right entry of $\xi$.


2. Solve the ordinary differential equation $d\Phi = \Phi \xi$.

3. Iwasawa split $\Phi$ into $\Phi = F B$.
Then $F$ is an extended frame for some CMC 
immersion.

4. Apply one of three Sym-Bobenko formulas described below to obtain
a CMC immersion into $\bbR^3$, $\bbS^3$ or $\bbH^3$.

\Subsection{The Sym-Bobenko formulas}
\label{sec:sym}

The final step in the extended Weierstrass representation
is a Sym-Bobenko formula, which computes the immersion
into $\bbR^3$, $\bbS^3$ or $\bbH^3$ from its extended frame.

1. CMC immersions into $\bbR^3$:
The Sym-Bobenko formula~\cite{Bobenko:cmc}
\begin{equation}
\label{eq:Sym-Bob}
- 2 i \lambda H^{-1} (\tfrac{d}{d\lambda} F) F^{-1}
\end{equation}
gives
a conformal CMC $H \neq 0$ immersion into $\bbR^3$ for each fixed
$\lambda_0 \in \bbS^1$. Formula~\eqref{eq:Sym-Bob} gives an
immersion into the Lie algebra $\matsu{2}{}$, which being a
real $3$-dimensional vector space can be identified with $\bbR^3$.

2. CMC immersions into $\bbS^3$:
For $\mu \in \bbS^1 \setminus \{ 1 \}$ and
each $\lambda_0 \in \bbS^1$, the Sym-Bobenko formula~\cite{Bobenko:cmc}
\begin{equation}
\label{eq:Sym-Bob-s3}
F_{\mu\lambda_0} F_{\lambda_0}^{-1}
\end{equation}
gives a conformal CMC $H = i (1+\mu)/(1-\mu)$ immersion into $\bbS^3$.
Here $H$ can take any real value, including $0$.
Formula~\eqref{eq:Sym-Bob-s3} gives an immersion into the Lie group
$\matSU{2}{}$, which
we are identifying with the unit sphere $\bbS^3 \in \bbR^4$.

3. CMC immersions into $\bbH^3$: For
$s \in (0,1)$ and any
$\lambda \in \bbS^1$, set $\lambda_0=s \lambda$.  Then the Sym-Bobenko
formula~\cite{Bobenko:cmc}
\begin{equation}
\label{eq:Sym-Bob-h3}
F_{\lambda_0} \transpose{\overline{F_{\lambda_0}}}
\end{equation}
gives a conformal CMC $H = (1+s^2)/(1-s^2) > 1$ immersion into $\bbH^3$
for each fixed $\lambda \in \bbS^1$.
Formula~\eqref{eq:Sym-Bob-h3} gives an immersion into the determinant
$1$ Hermitian matrices, which we are identifying with $\bbH^3$.

We choose the following normalizations for the Sym-Bobenko formulas:
\begin{equation}\label{eq:sym-normalize}
\begin{split}
&\text{
for $\bbR^3$: $\lambda_0=1$\spacecomma}\\
&\text{
for $\bbS^3$: $\lambda_0\in\bbS^{1}\setminus\{\pm 1\}$,
and $\mu = \lambda_0^{-2}$\spacecomma}\\
&\text{
for $\bbH^3$: $\lambda_0\in (-1,0) \cup (0,1)$\spaceperiod}
\end{split}
\end{equation}

\Subsection{Monodromy}
The primary result in~\cite{DPW} is that every CMC immersion into
$\bbR^3$ can be obtained via the extended Weierstrass representation,
and this is true for the cases of $\bbS^3$ and $\bbH^3$ as
well~\cite{SKKR:spaceforms} (with the restrictions $H \neq 0$ for
$\bbR^3$ and $|H|>1$ for $\bbH^3$).  This method can also be applied
to constructing non-simply-connected CMC immersions, and it is shown
in~\cite{DPW} that even when $\Sigma$ is a non-simply-connected
open non-compact Riemann surface, one can still always choose $\xi$ to
be well-defined on $\Sigma$, as long as the resulting CMC immersion is
well-defined on $\Sigma$.

However, in this case,
closing conditions must be satisfied in order for the resulting CMC
immersion to be well-defined on $\Sigma$.  Considering a deck
transformation $\tau$ of $\Sigma$ associated to some loop $\gamma$ in
$\Sigma$, let us suppose that we can choose the solution $\Phi$ so that
$\Phi \circ \tau = M_\gamma \Phi$ with $M_\gamma \in \LoopuSL{}$.
Then $M_\gamma$ is the monodromy of $\Phi$ about $\gamma$, and $M_\gamma$
is independent of $z$.  Because $M_\gamma \in \LoopuSL{}$, we also have
$F \circ \tau = M_\gamma F$.  Then the immersion obtained from the
Sym-Bobenko formula at $\lambda_0$ for $\bbR^3$ will be invariant
about $\gamma$ if $M_\gamma|_{\lambda_0} = \pm \id$ and
$\tfrac{d}{d\lambda} M_\gamma|_{\lambda_0} = 0$.  There are similar
conditions for the cases of $\bbS^3$ and $\bbH^3$. This gives the
following sufficient conditions for the resulting CMC immersion to be
well-defined on $\Sigma$:
\begin{align}
\label{eq:closing1}
&\text{for $\bbR^3$:
$M_\gamma \in \LoopuSL{}$,
$M_\gamma|_{\lambda_0} = \pm \id$,
$\tfrac{d}{d\lambda} M_\gamma|_{\lambda_0} = 0$\spacecomma}\\
\label{eq:closing2}
&\text{for $\bbS^3$:
$M_\gamma \in \LoopuSL{}$,
$M_\gamma|_{\mu \lambda_0} = M_\gamma|_{\lambda_0} = \pm \id$\spacecomma}\\
\label{eq:closing3}
&\text{for $\bbH^3$:
$M_\gamma \in \LoopuSL{}$,
$M_\gamma|_{\lambda_0} = \pm \id$}
\end{align}
for all loops $\gamma$ in $\Sigma$.  These are the
conditions we will show are
satisfied, to prove the existence of CMC trinoids and symmetric $n$-noids, by
making an appropriate choice of solution $\Phi$ of $d\Phi = \Phi \xi$.

\begin{longversion}
As the following easy lemma shows,
simultaneously unitarizing a finitely-generated group is equivalent
to unitarizing a set of generators of the group.

\begin{lemma}
\label{thm:group-uni}
Let $\calG\subset\matSL{2}{\bbC}$ be a finitely-generated subgroup.
If $C\in\matSL{2}{\bbC}$ unitarizes a set of generators of $\calG$,
then $C$ unitarizes every element of $\calG$.
\end{lemma}

\begin{proof}
1. If $C$ unitarizes any $X\in\matSL{2}{\bbC}$, then $C$ unitarizes $X^{-1}$.

2. If $C$ unitarizes any $X,\,Y\in\matSL{2}{\bbC}$, then $C$ unitarizes
$XY$.

3. Any element $P\in\langle M_1,\dots,M_n\rangle$ a product of the
form $P=P_1\dots P_n$, where the $P_k\in\{M_1^{\pm 1},\dots,M_n^{\pm
1}\}$.  By 1 and 2, $C$ unitarizes $P$.
\end{proof}
\end{longversion}

In the remainder of the paper we construct families of
CMC trinoids and symmetric $n$-noids.
For each family, the hypotheses of the Unitarization
Theorem~\ref{thm:uni-2by2} are shown to hold under a suitable set of
constraints on the end weights. The unitarization theorem then produces a
dressing which closes the ends.

\typeout{===nnoid=============}
\Section{Constructing $n$-noids}

\Subsection{The \textit{n}-noid potential}
\label{sec:nnoid-potential}

We define a class of potentials whose local
monodromies have the same eigenvalues as those of a Delaunay surface.

\begin{definition}
\label{def:potential}
Let $\Sigma=\bbP^1$ be the Riemann sphere
with the standard coordinate $z\in\bbC\cup\{\infty\}$.
Let $\lambda_0\in(\bbR\cup\bbS^1)\setminus\{0,-1\}$ be as
in~\eqref{eq:sym-normalize}, and let
\begin{equation}
\label{eq:h}
h(\lambda)
=
\fourth\lambda^{-1}(\lambda-\lambda_0)(\lambda-\lambda_0^{-1}) \spaceperiod
\end{equation}
Let $Q$ be a meromorphic quadratic differential on $\Sigma$
all of whose poles are double poles with real quadratic residues.
Assume that for each pole of $Q$, with quadratic residue $w/4$,
the function $1+ w h$ is non-negative on $\bbS^1$.
An \emph{$n$-noid potential} is an extended Weierstrass potential
of the form
\begin{equation*}
\label{eq:xi}
\xi=\begin{pmatrix}
0 & \lambda^{-1}dz \\ \lambda h(\lambda) Q/dz & 0
\end{pmatrix} \spaceperiod
\end{equation*}
\end{definition}

Let $p$ be a double pole of $Q$ with quadratic residue
$w/4\in\bbR\setminus\{0\}$.
Choosing a basepoint $z_0\in\Sigma$, let
$\gamma_p$ be a curve based at $z_0$ which winds
once around $p$ and does not wind
around any other poles of $Q$.
Let $M_p$
be the monodromy about $\gamma_p$
of the solution $\Phi=\Phi(z,\,\lambda)$
to the equation
$d\Phi=\Phi\xi$, $\Phi(z_0,\,\lambda)=\id$
along $\gamma_p$.

Multiplying $\Phi$ on the right by an analytic matrix $g=g(\lambda,z)$
does not change the resulting CMC immersion if $g = g_1 \cdot g_2$,
where $g_1$ is a $\lambda$-independent diagonal matrix
and $g_2 \in \LooppSL{}$.  We then call $\Phi g$
a \emph{gauge}
of $\Phi$.  This gauge
will change $\xi$ to
\begin{equation*}
\xi.g = g^{-1} \xi g + g^{-1} dg \spaceperiod
\end{equation*}

If $Q$ is holomorphic at $z=\infty$, then $\xi$ has a pole there.
The following lemma shows that this is an artifact of
our choice of potential,
not a feature of the monodromy representation or induced CMC immersion.
This lemma will be used in Section~\ref{sec:starnoid}.

\begin{lemma}\label{lem:gauge-inf}
Let $\xi$ be an $n$-noid potential as in Definition~\ref{def:potential} with 
$\Sigma=\bbP^1$.
Suppose $Q$ is holomorphic at $z=\infty\in\bbP^1$, and let $M_\infty$ be
a local monodromy at $\infty$.
Then 
$M_\infty=\id$, and 
$\infty$ is a smooth finite point of the CMC immersions
induced by the extended Weierstrass representation obtained from $\xi$.
\end{lemma}

\begin{proof}
Applying the gauge 
\begin{equation*}
g=\begin{pmatrix} z & 0 \\ -\lambda & z^{-1} \end{pmatrix} \spacecomma
\end{equation*}
the result follows from the fact that
$\gauge{\xi}{g}$ is holomorphic at $\infty$.
\end{proof}

\Subsection{Delaunay monodromy}

We will need several facts about the
$n$-noid monodromy defined in Section~\ref{sec:nnoid-potential}.
The first lemma computes the eigevalues of the monodromy and
proves the the latter half of the closing
conditions~\ref{eq:closing1}--\ref{eq:closing3}.

\begin{proposition}
\label{lem:monodromy-properties}
Let $M_p$ be a monodromy arising from an $n$-noid potential
as in Section~\ref{sec:nnoid-potential}. Then
\begin{enumerate}
\item
The eigenvalues of $M_p$ are $\exp(\pm 2\pi i \rho_w)$, where
\begin{equation}
\label{eq:rho}
\rho_w(\lambda)=\half-\half\sqrt{1+w h(\lambda)} \spaceperiod
\end{equation}
\item
With $\lambda_0$ as in~\eqref{eq:sym-normalize},
\begin{equation}
\label{eq:simple-closing}
M_p(\lambda_0^{\pm 1})=\id
\quad and\quad
\text{if $\lambda_0=1$, 
then $\left.\tfrac{d}{d\lambda}\right|_{\lambda_0}M_p=0$\spaceperiod}
\end{equation}
\end{enumerate}
\end{proposition}

\begin{proof}
The eigenvalues of $M_p$ can be computed using the theory
of regular singularities \cite{SKKR:spaceforms,DorW:tri}.

The first part of~\eqref{eq:simple-closing},
$M_p(\lambda_0^{\pm 1})=\id$, can be computed directly 
as the monodromy associated to $\xi(\lambda_0^{\pm 1})$.
To show the second part of~\eqref{eq:simple-closing}, assume $\lambda_0=1$.
Define the parameter $\theta$ by $\lambda =
e^{i\theta}$
and
let $L=\diag(e^{i\theta},\,e^{-i\theta})$.
Then $\xi$ has the symmetry
$\xi(-\theta)=L(\theta)\xi(\theta)L^{-1}(\theta)$,
from which it follows that
$M_p(-\theta)=L(\theta)M_p(\theta)L^{-1}(\theta)$.
Then
\begin{equation*}
0 = \left.\tfrac{d}{d\theta}\right|_{\theta=0} 
\left(M_p(-\theta)L(\theta)-L(\theta)M_p(\theta)\right)
= -2 \left. (\tfrac{d}{d\theta} M_p) \right|_{\theta=0} \; . 
\end{equation*}
The result
$\left.\tfrac{d}{d\lambda}\right|_{\lambda_0}M_p=0$ follows.
\end{proof}

\begin{longversion}

\begin{proof}
Let $M$ be the monodromy with basepoint $b$ of a solution
$\Phi_1$ to the equation
\[
d\Phi_1 = \Phi_1 \xi,
\quad
\Phi_1(b)=\id.
\]

Let $z_0$ be a pole of $Q$.
Assume without loss of generality by a change of coordinates
that $z_0=0$.
Then the series expansion of $\xi$ in $z$ at $0$ is
\[
\xi/dz = \begin{pmatrix} 0 & \lambda^{-1} \\
 \fourth\lambda hw & 0\end{pmatrix}z^{-2}
+\Order(x^{-1}).
\]
Define the gauge $g=\diag(z^{1/2},\,z^{-1/2})$.
Then the series expansion of $\xi$ in $z$ at $0$ is
\[
\xi/dz = A z^{-1} + \Order(z^0),
\]
where
\[
A = \begin{pmatrix}\half & \lambda^{-1} \\
 \fourth\lambda hw & -\half\end{pmatrix}
\]
Then, $\Phi_2 = g^{-1}(b) \Phi_1$ is a solution to the equation
\[
d\Phi_2 = \Phi_2 (\gauge{\xi}{g}),
\quad
\Phi_2(b) = \id.
\]

The eigenvalues of $A$ are $\half-\rho_w$,
so away from the set of isolated points at which $\rho_w\in\half\bbZ$,
we have by the theory of regular singularities
a holomorphic map $P$ such that
\[
\gauge{\xi}{P} = A.
\]
Then $\Phi_3 = \exp(A\log(z)) = P^{-1}(b)\Phi_2 P$
is a solution to the equation
\[
d\Phi_3 = \Phi_3 (\gauge{\xi}{g}),
\quad
\Phi_3(b) = \id.
\]
But the monodromy of $\Phi_3$ is $\exp(2\pi i A)$,
with eigenvalues $\exp(2\pi i(\half-\rho_w))$,
so the monodromy of $\Phi_1$
has eigenvalues $\exp(2\pi i\rho_w)$ away from 
the set of isolated points at which $\rho_w\in\half\bbZ$.

Since $\exp(2\pi i\rho_w)$ is a local analytic function,
and $M$ is analytic, it follows that the eigenvalues of
$M$ are $\exp(2\pi i\rho_w)$ at these isolated points.
\end{proof}
\end{longversion}

\Subsection{Local diagonalizability}

We show that subject to a bound on the end weights,
the $n$-noid monodromies satisfy
the local diagonalizability condition (ii)
of the Unitarization Theorem~\ref{thm:uni-2by2}.

\begin{lemma}
\label{lem:local-diagonalizability-lemma}
Let $\calC$ be an analytic curve and
$M\in\Holo{\calC}{\matM{2\times 2}{\bbC}}$
with local analytic eigenvalues $\mu_1,\mu_2\in\Holo{\calC}{\bbC}$
at $p\in \calC$. Then
\begin{enumerate}
\item
$\ord_p(\mu_1-\mu_2) \ge \ord_p(M-\mu_1\id)=\ord_p(M-\mu_2\id).$
\item
Assume $M$ is not identically a scalar multiple of $\id$.
Then the eigenlines of $M$ are non-coincident at $p$
if and only if
$\ord_p(\mu_1-\mu_2) = \ord_p(M-\mu_1 \id).$
\end{enumerate}
\end{lemma}

\begin{proof}
For
$M=\bigl(\begin{smallmatrix}a & b\\c & d\end{smallmatrix}\bigr)$,
define
$\adjoint{(M)}=\bigl(\begin{smallmatrix}d & -b\\-c & a\end{smallmatrix}\bigr)$.
To prove (i),
since $M+\adjoint{(M)} = (\mu_1+\mu_2)\id$, then
\begin{equation*}
M-\mu_1\id = -(\adjoint{(M)}-\mu_2\id) = -\adjoint{(M-\mu_2\id)} \spaceperiod
\end{equation*}
Hence
\begin{equation*}
\ord_p(M-\mu_1\id) = \ord_p(-\adjoint{(M-\mu_2\id)}) =
\ord_p(M-\mu_2\id).
\end{equation*}
Then, using $(\mu_1-\mu_2)\id = (M-\mu_2\id)-(M-\mu_1\id)$, we have
\begin{equation*}
\ord_p(\mu_1-\mu_2)\ge \min( \ord_p(M-\mu_1\id),\,\ord_p(M-\mu_2\id) )
= \ord_p(M-\mu_1\id) \spaceperiod
\end{equation*}

To prove (ii),
let $t$ be a local coordinate at $p$ on $\calC$ such that $t=0$ at $p$.
By (i), we can define
$n=\ord_p(M-\mu_1\id) = \ord_p(M-\mu_2\id)$.
Write
\begin{equation*}
M-\mu_k\id = A_k t^n + \Order(t^{n+1}) \spacecomma
\quad k\in\{1,\,2\}
\end{equation*}
with $A_1\ne 0$ and $A_2\ne 0$.
For $k\in\{1,\,2\}$,
the eigenline map $\calC\to\bbP^1$ corresponding to $\mu_k$ can
be written locally as $[v_k]$, where $v_k=a_k + \Order(t)$ for some
$a_k\in\bbC^2\setminus\{0\}$.
Then $(M-\mu_k\id)v_k=0$ implies $a_k\in\ker A_k$.
Then
\begin{equation*}
(\mu_1-\mu_2)\id = (M-\mu_2\id)-(M-\mu_1\id) = (A_2-A_1)t^n+\Order(t^{n+1})
\spacecomma
\end{equation*}
so $\ord_p(\mu_1-\mu_2)>n$ if and only if $A_1=A_2$.
If $A_1=A_2$, then $a_1\in\ker A_1$ and $a_2\in\ker A_1$,
and $a_1\ne 0$, $a_2\ne 0$, $A_1\ne 0$ imply $[a_1]=[a_2]$.
Conversely, if $[a_1]=[a_2]$, then
$a_1\in\ker A_1$ and $a_1\in\ker A_2$, so
$a_1\in\ker(A_1-A_2)$.
Since $A_1-A_2$ is a scalar multiple of $\id$,
then $A_1-A_2=0$.
\end{proof}

\begin{lemma}
\label{lem:noid-local-diag}
Let $M=M_p \in \LoopSL{}$ be an $n$-noid monodromy at $p$ as above.
Let $\rho=\rho_w$ be as in~\eqref{eq:rho}
and assume $\abs{\rho}<\half$ on $\bbS^1$.
Then $M$ is locally diagonalizable at each point of $\bbS^1$.
\end{lemma}

\begin{proof}
Let $\mu$ be an eigenvalue of $M$.
Then $\mu$ is locally analytic on $\bbS^1$ because
$\half\tr M\in[-1,\,1]$ on $\bbS^1$.
Let $\lambda_0$ be as in~\eqref{eq:sym-normalize}.
Because $\abs{\rho}<\half$, $\mu$ is never $-1$ on $\bbS^1$,
and $\mu$ is $+1$ on $\bbS^1$ only at $\lambda_0^{\pm 1}$.
Define $n=n_{\lambda_0}:\bbS^1\to\{0,\,1,\,2\}$ by
\begin{align*}
&\text{$n_{\lambda_0}(p) = 0$ if $p\in\bbS^1\setminus\{\lambda_0^{\pm 1}\}$\spacecomma}\\
&\text{$n_{\lambda_0}(\lambda_0^{\pm 1})=1$ if $\lambda_0\in\bbS^1\setminus\{1\}$\spacecomma}\\
&\text{$n_{\lambda_0}(\lambda_0)=2$ if $\lambda_0=1$\spaceperiod}
\end{align*}
Then for all $p\in\bbS^1$, we have
$\ord_p(\mu-1) = n(p) = \ord_p(\mu-\mu^{-1})$,
and by~\eqref{eq:simple-closing},
$\ord_p(M-\id)\ge n(p)$,
Then using $M-\mu\id = (M-\id) - (\mu - 1)\id$, we have
\begin{equation*}
\ord_p(M-\mu\id)
 \ge
 \min(\ord_p(M-\id),\,\ord_p(\mu-1))
 = n(p)
 = \ord_p(\mu-\mu^{-1})\spacecomma
\end{equation*}
and the result follows by Lemma~\ref{lem:local-diagonalizability-lemma}.
\end{proof}

\typeout{===nnoid3=============}
\Subsection{Unitarizing three loops whose product is $\id$}
\label{sec:nnoid3}

The following well-known proposition \cite{Goldman:top}
(see also~\cite{Biswas:pun,Umehara-Yamada})
gives a condition for simultaneous unitarizability of
three matrices whose product is $\id$ in terms
of their traces.

\begin{proposition}[\cite{Goldman:top}]
\label{lem:goldman-noid}
Let $M_1,\,M_2,\,M_3\in\matSL{2}{\bbC}$
and suppose $M_1M_2M_3=\id$.
For $k\in\{1,\,2,\,3\}$,
let $t_k = \half\tr M_k$ and suppose $t_k\in[-1,\,1]$.
Define
\begin{equation}
\label{eq:T-X}
T = 1-t_1^2-t_2^2-t_3^2+2t_1t_2t_3 \spaceperiod
\end{equation}
Then
the $M_1,\,M_2,\,M_3$ are
reducible if and only if $T=0$,
and are
irreducible and simultaneously unitarizable
if and only if $T>0$.
\end{proposition}

If $M_1,\,M_2,\,M_3:\calC\to\matSL{2}{\bbC}$ are analytic maps on
an analytic curve $\calC$,
then their infinitesimal irreducibility at a zero of $T$
can in some cases be computed by the following technique.

\begin{lemma}
\label{lem:det-T-X}
Let $\calC$ be an analytic curve and $p\in\calC$.
Let $M_1,\,M_2,\,M_3\in\Holo{\calC}{\matSL{2}{\bbC}}$
satisfy $M_1M_2M_3=\id$.
Let $T$ be 
as in ~\eqref{eq:T-X}
with $t_k =\half\trace M_k$, $k\in\{1,\,2,\,3\}$.
Then for any pair $j,\,k$ in $\{1,\,2,\,3\}$, 
if $\ord_p [M_j,\,M_k] \ge n \ge 0$,
\begin{equation*}
\det([M_j,\,M_k]^{(n)})(p) = \tfrac{4}{b_{2n,n}} T^{(2n)}(p) \spacecomma
\end{equation*}
where the superscript $(n)$ denotes differentiation $n$ times
with respect to a coordinate at $p$, and $b_{r,s}$ denotes the
binomial coefficient.
\end{lemma}

\begin{proof}
We first show that for any analytic map 
$X: \calC \to\matM{2\times 2}{\bbC}$ with $\trace X \equiv 0$,
if $\ord_p X \ge n \ge 0$,  
then
\begin{equation}
\label{eq:det-T1-X}
\det( X^{(n)}(p) ) = \tfrac{1}{b_{2n,n}}(\det X)^{(2n)}(p) \spaceperiod
\end{equation}
Differentiating the Cayley-Hamilton equation $(\det X) \id = -X^2$, we have
\begin{equation*}
(\det X)^{(2n)} \cdot \id
 =
 -\textstyle\sum_{k=0}^{2n}b_{2n,k}X^{(k)}X^{(2n-k)} \spaceperiod
\end{equation*}
As $X^{(k)}(p)=0$ for $k < n$, all terms in the sum are zero except
possibly 
when $k=n$. This yields
\begin{equation*}
(\det X)^{(2n)}(p) \cdot \id = -b_{2n,n}(X^{(n)}(p))^2 =
b_{2n,n}\det( X^{(n)}(p) ) \cdot \id \spacecomma
\end{equation*}
proving~\eqref{eq:det-T1-X}.

A calculation shows
\begin{equation}
\label{eq:det-T2-X}
\det([M_j,M_k])=4T \spaceperiod
\end{equation}
The result follows directly
from~\eqref{eq:det-T2-X} 
and~\eqref{eq:det-T1-X} with $X=[M_j,\,M_k]$.
\end{proof}

\typeout{===trinoid=============}

\Section{Trinoids}
\label{sec:3noid}

In~\cite{SKKR:spaceforms} a three-parameter family of constant mean curvature
trinoids was constructed for each mean curvature $H$ in each of the
space forms $\bbR^3$, $\bbS^3$ and $\bbH^3$,
using $r$-Iwasawa decomposition for $r<1$.
Here we show, employing the $1$-unitarization Theorem~\ref{thm:uni-2by2},
that these immersions can be constructed with
less machinery, using only the $1$-Iwasawa decomposition.

\Subsection{Trinoid potentials}
\label{sec:3noid-pot}

\begin{definition}
[Trinoid potentials]
\label{def:3noid}
Let $\Sigma=\bbP^{1}\setminus\{0,1,\infty\}$ be the thrice-punctured
Riemann sphere.
The family of \emph{trinoid potentials} $\xi$ on $\Sigma$,
parametrized by
$\lambda_0$ and
$w_0,w_1,w_\infty \in\bbR\setminus\{0\}$,
is given by $\xi$ in Definition~\ref{def:potential} with
\begin{equation*}
Q=\frac{w_\infty z^2+(w_1-w_0-w_\infty)z+w_0}{4z^2(z-1)^2}dz^2 \spaceperiod
\end{equation*}
\end{definition}

$Q$ is the unique quadratic differential with
double poles at $\{0,\,1,\,\infty\}$
(the ends of the surface)
with respective quadratic residues
$w_0/4,\,w_1/4,\,w_\infty/4$,
and no other poles.

A set of generators of the monodromy representation of a
trinoid potential $\xi$ is defined as follows.
Choose a basepoint $z_0\in\bbP^{1}\setminus\{0,1,\infty\}$.
For $k\in\{0,\,1,\infty\}$,
a set of closed curves $\gamma_k$
based at $z_0$ can be chosen which wind respectively
around $k\in\bbP^{1}$ once and not
around any other point in $\{0,\,1,\infty\}$, satisfying
$\gamma_0\gamma_1\gamma_\infty=\id$.
Define $M_k:\bbC^\ast\to\matSL{2}{\bbC}$ as the monodromy
of the solution $\Phi = \Phi(z,\,\lambda)$ to the equation
$d\Phi=\Phi\xi$, $\Phi(z_0,\,\lambda)=\id$ along $\gamma_k$.
Then by the choice of $\gamma_0$, $\gamma_1$, $\gamma_\infty$,
we have $M_0M_1M_\infty=\id$.

\Subsection{Pointwise unitarizability}
\label{sec:3noid-uni}

A key step from~\cite{SKKR:spaceforms} in the trinoid construction
is showing, with a suitable set of inequalities, that the
monodromy representation is pointwise unitarizable on $\bbS^1$.
The following lemma is a restatement of the required
lemma in~\cite{SKKR:spaceforms}
with inequalities replaced by strict inequalities.

\begin{lemma}[\cite{SKKR:spaceforms}]
\label{lem:3noid-uni}
Let $\xi$ be a trinoid potential
parametrized by $\lambda_0,w_0,w_1,w_\infty$,
and let $\{M_0,M_1,M_\infty\}$ be the generators of the monodromy
representation for $\xi$ as described above.
For $k\in\{0,1,\infty\}$ define
$\rho_k = \rho_{w_k}$ as in~\eqref{eq:rho} and
$n_k = \rho_{w_k}(-1)$ and $m_k = \rho_{w_k}(1)$.
Suppose the following inequalities hold
for every permutation $(i,j,k)$ of $(0,1,\infty)$:
\begin{align}
\label{eq:3noid-inequality1}
&\abs{n_0}+\abs{n_1}+\abs{n_\infty} < 1
\text{ and }
\abs{n_i} < \abs{n_j}+\abs{n_k}
\text{ for all space forms,}\\
\label{eq:3noid-inequality2}
&\abs{m_0}+\abs{m_1}+\abs{m_\infty} < 1
\text{ and }
\abs{m_i} < \abs{m_j}+\abs{m_k}
\text{ for $\bbS^3$ and $\bbH^3$\spacecomma}\\
\label{eq:3noid-inequality3}
&\abs{w_i} < \abs{w_j}+\abs{w_k}
\text{ for $\bbR^3$\spaceperiod}
\end{align}
Then the monodromy representation for $\xi$
is pointwise unitarizable on $\bbS^1$, and is
irreducible on $\bbS^1\setminus\{\lambda_0^{\pm 1}\}$.
Moreover, $\abs{\rho_k}<\half$ on $\bbS^1$.
\end{lemma}

\Subsection{Infinitesimal irreducibility}
\label{sec:3noid-irred}

\begin{lemma}
\label{lem:3noid-irred}
Let $\xi$ be a trinoid potential
parametrized by $\lambda_0$, $w_0,\,w_1,\,w_\infty$.
Let $M_0,M_1,M_\infty$ be the monodromies for $\xi$
as in Section~\ref{sec:3noid-uni}.
Suppose condition~\eqref{eq:3noid-inequality3} holds.
Then $M_0,M_1,M_\infty$ are pairwise
infinitesimally irreducible at $\{\lambda_0^{\pm 1}\}$.
\end{lemma}

\begin{proof}
Let $'$ denote differentiation with respect to $\lambda$,
and let the superscript $(k)$ denote differentiation $k$ times with
respect to $\lambda$. Let $T$ be as in~\eqref{eq:T-X}.
Choose distinct $j,k\in\{0,1,\infty\}$.

Let $\chi\in\bbR$ be defined by
\begin{equation}
\label{eq:chi}
\chi=(w_0+w_1+w_\infty)(-w_0+w_1+w_\infty)(w_0-w_1+w_\infty)(w_0+w_1-w_\infty)
\end{equation}
Note that $\chi=0$ if and only if
$\abs{w_i}=\abs{w_j}+\abs{w_k}$ for some permutation $(i,\,j,\,k)$ of
$(0,\,1,\,\infty)$.
Hence by condition~\eqref{eq:3noid-inequality3}, $\chi\ne 0$.

First take the case $\lambda_0\ne 1$.
A calculation using $M_r(\lambda_0^{\pm 1})=\id$, $r\in\{j,\,k\}$, shows
$\ord_{\lambda_0^{\pm 1}}[M_j,\,M_k]\ge 2$.
By Lemma~\ref{lem:det-T-X} and a calculation,
\begin{equation*}
\begin{split}
\det([M_j,M_k]^{(2)}(\lambda_0^{\pm 1}))
&=
\tfrac{4}{b_{4,2}} T^{(4)}(\lambda_0^{\pm 1})
=
\tfrac{4}{b_{4,2}} 3\cdot 2^{-11}\pi^4 (1-\lambda_0^{\mp 2})^4 \chi \spaceperiod
\end{split}
\end{equation*}
Thus since $\chi\ne 0$,
then $\det([M_j,M_k]^{(2)}(\lambda_0^{\pm 1}))\ne 0$
and $M_j$ and $M_k$ are infinitesimally irreducible at
$\lambda_0^{\pm 1}$
by Definition~\ref{def:infinitesimally-irreducible}.

For the case $\lambda_0=1$,
a calculation using
 $M_r(1)=\id$,
 $M_r'(1)=0$,
 $r\in\{j,\,k\}$, shows
$\ord_{1}[M_j,\,M_k]\ge 4$.
By Lemma~\ref{lem:det-T-X} and a calculation,
\begin{equation*}
\begin{split}
\det([M_j,M_k]^{(4)}(1))
&=
\tfrac{4}{b_{8,4}} T^{(8)}(1)
=
\tfrac{4}{b_{8,4}}315\cdot 2^{-7}\pi^4 \chi \spaceperiod
\end{split}
\end{equation*}
Thus since $\chi\ne 0$,
then $\det([M_j,M_k]^{(4)}(1))\ne 0$
and $M_j$ and $M_k$ are infinitesimally irreducible at $1$
by Definition~\ref{def:infinitesimally-irreducible}.
\end{proof}

\Subsection{Constructing trinoids}
\label{sec:3noid-construct}

The results of Sections~\ref{sec:3noid-pot}--\ref{sec:3noid-irred}
are now brought together to construct a family of trinoids.

\begin{theorem}[Trinoids]
\label{thm:3noid}
Let $\xi$ be a trinoid potential
on $\Sigma=\bbP^1\setminus\{0,\,1,\,\infty\}$
satisfying the inequalities
\eqref{eq:3noid-inequality1}--\eqref{eq:3noid-inequality3}.
Then there exists a solution $\Psi$ of the equation $d\Psi=\Psi \xi$
such that $\Psi$ induces a CMC $H$ immersion of $\Sigma$
into the appropriate space form $\bbR^3$ or $\bbS^3$ or $\bbH^3$
via the extended Weierstrass representation,
where the mean curvature $H$ is subject to the restrictions in
Section~\ref{sec:sym}.
\end{theorem}

\begin{proof}
Let $z_0\in\Sigma$ be a basepoint, and let $\Phi$ be the solution
to $d\Phi=\Phi\xi$, $\Phi(z_0,\,\lambda)=\id$.
We show that the hypotheses of Theorem~\ref{thm:uni-2by2} hold
for the generators
$\{M_0,\,M_1,\,M_\infty\}$
of the trinoid monodromy representation.

By Lemma~\ref{lem:3noid-uni},
$M_0,\,M_1,\,M_\infty$ are pairwise irreducible on
$\bbS^1\setminus\{\lambda_0^{\pm 1}\}$ and hence no two identically
commute.
By the same lemma,
$M_0,\,M_1,\,M_\infty$ are pointwise unitarizable on $\bbS^1$,
condition (i) of Theorem~\ref{thm:uni-2by2}.

%
Lemma~\ref{lem:3noid-uni} provides the bound $\abs{\rho_k}<\half$,
so by Lemma~\ref{lem:noid-local-diag},
$M_0,\,M_1,\,M_\infty$ are each locally diagonalizable at
each point of $\bbS^1$,
condition (ii) of Theorem~\ref{thm:uni-2by2}.

Since
$M_0,M_1,M_\infty$ are irreducible on
$\bbS^1\setminus\{\lambda_0^{\pm 1}\}$,
they are pairwise infinitesimally irreducible there.
By Lemma~\ref{lem:3noid-irred},
$M_0,M_1,M_\infty$ are pairwise infinitesimally irreducible at
$\{\lambda_0^{\pm 1}\}$.
Therefore
$M_0,M_1,M_\infty$ are pairwise infinitesimally irreducible on $\bbS^1$,
condition (iii) of Theorem~\ref{thm:uni-2by2}.

Thus all conditions of Theorem~\ref{thm:uni-2by2} are
satisfied, so by that theorem there exists an analytic loop
$C \in \LoopSL{}$ which unitarizes the monodromy representation of $\Phi$.
In the case of the spaceform $\bbH^3$, $C$ may be singular at $\lambda_0$, so
let $C=C_u \cdot C_+$ be the 1-Iwasawa decomposition of $C$.
Then $C_+$ likewise unitarizes the monodromy representation of $\Phi$,
and is nonsingular at $\lambda_0$ for any spaceform.

Then $\Psi=C_+ \Phi$ satisfies the appropriate closing
condition~\eqref{eq:closing1}--\eqref{eq:closing3}
since condition~\eqref{eq:simple-closing}
is independent of conjugation by an analytic loop.
Hence the immersion induced
by $\Psi$ is well-defined on $\Sigma$.
\end{proof}

\begin{remark}
\label{prop:boundary_case}
In the case in which any of the
inequalities~\eqref{eq:3noid-inequality1}--\eqref{eq:3noid-inequality3}
becomes an equality,
the trinoid can be constructed via
the $r$-Unitarization Theorem~\ref{thm:runi}
(see~\cite{SKKR:spaceforms}).
\end{remark}

\typeout{===starnoid=============}

\Section{Symmetric $n$-noids}
\label{sec:starnoid}

The above ideas are now applied to the construction of CMC symmetric $n$-noids,
genus-zero surfaces similar to trinoids, but having $n$ ends
and full dihedral symmetry of order $n$. 

\Subsection{Symmetric $n$-noid potentials}
\label{sec:nnoid-pot}

\begin{definition}[Symmetric $n$-noid potentials]
\label{def:nnoid}
Let $n$ be an integer with $n\ge 3$.
Let $\Sigma$ be the $n$-punctured Riemann sphere
$\Sigma=\bbP^1\setminus\{z^n=1\}$.
Let $w\in\bbR\setminus\{0\}$.
The family of \emph{symmetric $n$-noid potentials} $\xi$,
parametrized
by $\lambda_0$ and $w$,
is given by $\xi$ in Definition~\ref{def:potential}
with
\begin{equation*}
Q=\frac{n^2 w z^{n-2}}{4(z^n-1)^2}dz^2 \spaceperiod
\end{equation*}
\end{definition}

$Q$ is chosen to have double poles at $\{z^n=1\}$
(the ends of the surface),
with quadratic residue $w/4$ at each pole,
and no other poles.
The choice of ends gives $Q$ a symmetry that will imply
the pointwise simultaneously unitarizability
of the monodromy group for $\xi$ on $\bbS^1$.

A set of generators of the monodromy representation of a
symmetric $n$-noid potential $\xi$ is defined as follows.
A set of closed curves
$\gamma_0,\dots,\gamma_{n-1},\gamma_\infty$ based 
at $0$ can be chosen which respectively wind around 
$e^{2 \pi i 0/n},\dots,e^{2\pi i (n-1)/n},\,\infty$ once
and not
around any other of these points,
satisfying $\gamma_0\dots\gamma_{n-1}\gamma_\infty=\id$.
Define $M_0,\dots,M_{n-1},\,M_\infty:\bbC^\ast\to\matSL{2}{\bbC}$
as the monodromies
of the solution $\Phi(z,\,\lambda)$ to the equation
$d\Phi=\Phi\xi$, $\Phi(0,\,\lambda)=\id$ along
$\gamma_0,\dots,\gamma_{n-1},\,\gamma_\infty$ respectively.
This
choice gives
$M_0\cdots M_{n-1}M_\infty=\id$.

\begin{lemma}
\label{thm:nnoid-symmetry}
Let $\xi$ be a symmetric $n$-noid potential.
Define the gauge
\begin{equation}
\label{eq:alpha-g}
g=\diag(\alpha^{1/2},\,\alpha^{-1/2}),\quad \alpha=e^{2\pi i/n} \spaceperiod
\end{equation}
Let $\tau:\bbP^1\to\bbP^1$ be the automorphism of $\bbP^1$
defined by $\tau(z)=\alpha z$.
Then
\begin{enumerate}
\item $\xi$ has the symmetry $\tau^\ast\xi = \gauge{\xi}{(g^{-1})}$.
\item
Let $\Phi=\Phi(z,\,\lambda)$ solve
$d\Phi = \Phi \xi$, $\Phi(0,\,\lambda)=\id$.
Then $\Phi$ has the symmetry $\tau^\ast\Phi = g \Phi g^{-1}$.
\item The monodromy matrices $M_0,\dots,M_{n-1}$ of $\Phi$ satisfy
$M_k = g^k M_0 g^{-k}$ for $k=0,\dots,n-1$.
\end{enumerate}
\end{lemma}

\begin{proof}
Showing (i) is a calculation.
By (i), every solution of $d\Phi = \Phi\xi$ has the symmetry
$\tau^\ast\Phi = A \Phi g^{-1}$ for some $A$.
Evaluating at $z=0$ yields $A=g$, implying (ii).
Symmetry (iii) follows.
\end{proof}

\Subsection{Pointwise unitarizability}
\label{sec:nnoid-uni}

The techniques of Section~\ref{sec:nnoid3} do not directly apply
to the generators $M_0,\dots,M_{n-1}$ of 
the symmetric $n$-noid monodromy
representation, but rather to the triple $M_0$, $g$, $(M_0g)^{-1}$,
whose product is $\id$ and whose traces are computable.
The pointwise or loopwise unitarizability of this triple
implies the same for the symmetric $n$-noid monodromy representation.

The value of $w$ in the symmetric $n$-noid potential
determines the weight of the Delaunay potential to which the
symmetric $n$-noid potential is
asymptotic~\cite{KRS:asymptotics,SKKR:spaceforms,KorKMS,KorKS}.
So the condition~\eqref{eq:nnoid-inequality} in the lemma below
amounts to a restriction on the weight of the ends.

\begin{lemma}
\label{lem:nnoid-uni}
Let $\xi$ be an $n$-noid potential with parameters $w$ and $\lambda_0$.
With $h$ as in~\eqref{eq:h} and $\rho_w$ as in~\eqref{eq:rho},
let $M_0$ be the symmetric $n$-noid monodromy defined
above and let $g$ be 
as in~\eqref{eq:alpha-g}.
Suppose
\begin{align}
\label{eq:nnoid-inequality}
\abs{\rho(1)} < \tfrac{1}{n}
\quad\text{and}\quad
\abs{\rho(-1)} < \tfrac{1}{n} \spaceperiod
\end{align}
Then $g$ and $M_0$
are pointwise simultaneously unitarizable on $\bbS^1$,
and are irreducible on $\bbS^1\setminus\{\lambda_0^{\pm 1}\}$.
\end{lemma}

\begin{proof}
\begin{longversion}
Since $\rho$ attains it minimum and maximum on $\bbS^1$ at $\lambda=1$
and $\lambda=-1$ in some order,
condition~\eqref{eq:nnoid-inequality} is equivalent
to the condition that $\rho$ takes values in $(-1/n,\,1/n)$ on $\bbS^1$.
\par
\end{longversion}
By Lemma~\ref{lem:gauge-inf} we have $M_0\cdots M_{n-1}=\id$.
A calculation using
Lemma~\ref{thm:nnoid-symmetry}(iii) implies
\begin{equation}
\label{eq:Mg}
(M_0g)^n = -\id \; . 
\end{equation}
It follows that the eigenvalues of $M_0g$ are
$n$'th roots of $-1$, and are hence constant.
With $\alpha$ as in~\eqref{eq:alpha-g}, using $M_0(\lambda_0)=\id$ we get
\begin{equation}
\label{eq:trace-Mg}
\text{the eigenvalues of $M_0g$ are $\alpha^{\pm 1/2}$\spaceperiod}
\end{equation}

Now consider the triple $(M_0,\,g,\,(M_0g)^{-1})$.
Their product is $\id$, 
and
\begin{equation*}
t:=\half\tr M_0=\cos(2\pi\rho_w) \spacecomma
\quad
\half\tr g= \half\tr ((M_0g)^{-1})=(\alpha^{1/2}+\alpha^{-1/2})/2 \spaceperiod
\end{equation*}
Hence with $T$ as in~\eqref{eq:T-X}, 
\begin{equation}
\label{eq:nnoid-T}
T(\lambda)= (1-t)(t-(\alpha+\alpha^{-1})/2) \spaceperiod
\end{equation}

Condition~\eqref{eq:nnoid-inequality} implies that
$\rho_w$ takes values in $(-1/n,\,1/n)$ on $\bbS^1$,
and $\rho_w$ is zero only at $\lambda_0^{\pm 1}$.
Hence $t$ takes values in $((\alpha+\alpha^{-1})/2,\,1]$
and the zero set of $T$ on $\bbS^1$ is $\{\lambda_0^{\pm 1}\}$.
Then by Proposition~\ref{lem:goldman-noid},
($M_0,\,g,\,(M_0 g)^{-1})$ are simultaneously unitarizable
and irreducible on $\bbS\setminus\{\lambda_0^{\pm 1}\}$.
We have by~\eqref{eq:simple-closing}
that $M_0=\id$, so $(M_0,\,g,\,(M_0g)^{-1})$ are simultaneously unitarizable
at $\{\lambda_0^{\pm 1}\}$, and hence on $\bbS^1$.
\end{proof}

\begin{longversion}
Proof in other words:
Using~\eqref{eq:nnoid-inequality},
the zero set of $T$ on $\bbS^1$ is $\{\lambda_0^{\pm 1}\}$.
Then by Proposition~\eqref{lem:goldman-noid},
$M_0$, $g$, $(M_0 g)^{-1}$ are simultaneously unitarizable
and irreducible on $\bbS\setminus\{\lambda_0^{\pm 1}\}$.

At the exceptional point $\{\lambda_0^{\pm 1}\}$,
we have by~\eqref{eq:simple-closing}
that $M_0=\id$, so $g$ and $M$ are simultaneously unitarizable
at $\{\lambda_0^{\pm 1}\}$.
\end{longversion}

\Subsection{Infinitesimal irreducibility}
\label{sec:nnoid-irred}


\begin{lemma}
\label{lem:nnoid-irred}
Let $\xi$ be a symmetric $n$-noid potential.
Let $M_0$ be the symmetric $n$-noid monodromy defined
in Section~\ref{sec:nnoid-uni}
and $g$ as in~\eqref{eq:alpha-g}.
Then $g$ and $M_0$ are infinitesimally irreducible
at $\{\lambda_0^{\pm 1}\}$.
\end{lemma}

\begin{proof}
If $\lambda_0\ne 1$, then
$M_0(\lambda_0^{\pm 1})=\id$,
so $\ord_{\lambda_0^{\pm 1}}[g,\,M_0]\ge 1$.
By Lemma~\ref{lem:det-T-X} and a calculation using~\eqref{eq:nnoid-T}, 
taking derivatives with respect to the parameter $\theta$ for 
$\lambda = e^{i \theta}$, 
\begin{equation*}
\det([g,\,M_0]^{(1)}(\lambda_0^{\pm 1}))
 = \tfrac{4}{b_{2,1}} T^{(2)}(\lambda_0^{\pm 1})
 =
 \tfrac{1}{\alpha b_{2,1}} 2^{-5} \pi^2 (1-\alpha)^2 
 (\lambda_0-\lambda_0^{-1})^2 w^2 \spacecomma
\end{equation*}
and since $\lambda_0^2\ne 1$,
$\alpha\ne 1$
and $w\ne 0$, 
then $\det([g,\,M_0]^{(1)}(\lambda_0^{\pm 1}))\ne 0$
and 
$g$ and $M_0$ are infinitesimally irreducible
at $\lambda_0^{\pm 1}$
by the definition of infinitesimal irreducibility
in Section~\ref{sec:twobytwo}.

If $\lambda_0=1$, then
$M_0(1)=\id$ and $M_0^{(1)}(1)=0$,
so $\ord_1[g,\,M_0]\ge 2$.
By Lemma~\ref{lem:det-T-X} and a calculation using~\eqref{eq:nnoid-T},
\begin{equation*}
\det([g,\,M_0]^{(2)}(1))
=
\tfrac{4}{b_{4,2}} T^{(4)}(1)
=
-\tfrac{1}{\alpha b_{4,2}} 3\cdot 2^{-3} \pi^2 (1-\alpha)^2 w^2\spacecomma
\end{equation*}
so $\det([g,\,M_0]^{(2)}(1))\ne 0$
and 
$g$ and $M_0$ are infinitesimally irreducible at $1$.
\end{proof}

\Subsection{Constructing symmetric $n$-noids}
\label{sec:nnoid-construct}

The results of Sections~\ref{sec:nnoid-pot}--\ref{sec:nnoid-irred}
are now brought together to construct a family of symmetric $n$-noids.

\begin{theorem}
\label{thm:nnoid}
With $n\ge 3$,
let $\xi$ be a symmetric $n$-noid potential on
$\Sigma=\bbP^1\setminus\{z^n=1\}$ such that the
inequalities~\eqref{eq:nnoid-inequality} hold.
Then there exist solutions $\Psi$ of the equation $d\Psi = \Psi\xi$
such that the $\Psi$ induce
a 1-parameter family of
CMC $H$ immersions of $\Sigma$ into
each of the space forms $\bbR^3$, $\bbS^3$ and $\bbH^3$,
where $H\in\bbR$ is
subject to the restrictions of Section~\ref{sec:sym}.
\end{theorem}

\begin{proof}
%
Let $M_0$ be the monodromy described in Section~\ref{sec:nnoid-pot}
and let $g$ be as in~\ref{thm:nnoid-symmetry}.
Let $\Phi$ solve the equation $d\Phi=\Phi\xi$
with $\Phi(0,\,\lambda)=\id$.
We show that the hypotheses of Theorem~\ref{thm:uni-2by2} hold
for the
triple $\{M_0,\,g,\,{(M_0g)}^{-1}\}$.

By Lemma~\ref{lem:nnoid-uni}, $M_0$ and $g$ are irreducible on
$\bbS^1\setminus\{\lambda_0^\pm\}$, and hence $[M_0,\,g]\not\equiv 0$.
By the same lemma, $M_0$ and $g$ are pointwise
simultaneously unitarizable at every
point of $\bbS^1$,
condition (i) of Theorem~\ref{thm:uni-2by2}.
%
%
Condition (ii) of Theorem~\ref{thm:uni-2by2} holds because $g$ is diagonal.
Lemmas~\ref{lem:nnoid-uni} and~\ref{lem:nnoid-irred}
show 
that $M_0$ and $g$ are infinitesimally irreducible at every
point of $\bbS^1$,
condition (iii) of Theorem~\ref{thm:uni-2by2}.

Thus all conditions of Theorem~\ref{thm:uni-2by2} are
satisfied, so 
there exists an analytic loop
$C \in \LoopSL{}$ which simultaneously unitarizes
$M_0$ and $g$.
$C$ unitarizes the monodromy representation
described in Section~\ref{sec:nnoid-pot}, 
because it is contained in the group generated by $M_0$ and $g$.
In the case of the spaceform $\bbH^3$, $C$ may be singular at $\lambda_0$, so
let $C=C_u \cdot C_+$ be the 1-Iwasawa decomposition of $C$.
Then $C_+$ likewise unitarizes the monodromy representation.

Let $\Psi = C_+\Phi$.
Then the solution $\Psi$ to
$d\Phi=\Psi\xi$, $\Psi(0,\,\lambda)=C_+$ has unitary monodromy
satisfying the appropriate closing
condition~\eqref{eq:closing1}--\eqref{eq:closing3},
since condition~\eqref{eq:simple-closing} 
is independent of conjugation by an analytic loop.
Hence the CMC immersion induced
by $\Psi$ via the extended Weierstrass representation
into the appropriate space form
is well-defined on $\Sigma$.
Note that $z=\infty$ is a finite smooth point of the immersion,
by Lemma~\ref{lem:gauge-inf}, so the surface has $n$ ends.
For each spaceform and each choice of $n$,
this produces a one-parameter family of surfaces parametrized by $w$.
%
\end{proof}

\begin{longversion}
-----------------
Proof in other words.
\begin{proof}
Let $z_0\in\Sigma$ be a basepoint, and let $\Phi$ be the solution
to $d\Phi=\Phi\xi$, $\Phi(z_0,\,\lambda)=\id$.
We show that the hypotheses of Theorem~\ref{thm:uni-2by2} hold, applied
to the two matrix maps $g$ and $M_0$.
By Lemma~\ref{lem:nnoid-uni}, $[g,\,M_0]\not\equiv 0$.

The inequality~\eqref{eq:nnoid-inequality}
implies that $1+wh$ is non-zero on $\bbS^1$.
$M_0$ is locally diagonalizable by Lemma~\ref{lem:noid-local-diag},
condition (i) of Theorem~\ref{thm:uni-2by2}.

Lemma~\ref{lem:nnoid-uni} and~\ref{lem:nnoid-irred}
show that $M_0$ and $g$ are infinitesimally irreducible,
condition (ii) of Theorem~\ref{thm:uni-2by2}.

Lemma~\ref{lem:nnoid-uni} shows that $M_0$ and $g$ are pointwise
simultaneously unitarizable,
condition (iii) of Theorem~\ref{thm:uni-2by2}.

Thus all conditions of Theorem~\ref{thm:uni-2by2} are
satisfied, so by that theorem there exists an analytic loop
$C \in \LoopSL{}$ which unitarizes the monodromy representation of $\Phi$.
Let $C=C_u \cdot C_+$ be the 1-Iwasawa decomposition of $C$.
Then $C_+$ likewise unitarizes the monodromy representation of $\Phi$.

Then $\Psi=C_+ \Phi$ satisfies the appropriate closing
condition~\eqref{eq:closing1}--\eqref{eq:closing3}
since condition~\eqref{eq:simple-closing},
is independent of conjugation by an analytic loop.
Hence the immersion induced
by $C_+\Phi$ is is well-defined on $\Sigma$.
\end{proof}
\end{longversion}

\begin{remark}
The methods in this section can be extended to
a broader family of symmetric $n$-noids whose end axes are not coplanar
\cite{Schmitt:noid}.
\end{remark}


\bibliographystyle{amsplain}
\bibliography{unitarize}

\end{document}